\newif\ifarxivreport 
\arxivreporttrue 

\documentclass[10pt,a4paper,reqno]{amsart}

\newcommand{\boi}{\text{boi}}
\newcommand{\chiv}{\text{chi}}

\newcommand{\im}{\text{im}}

\newcommand{\cel}{\degree	\text{C}}

\newcommand{\ogam}{\overline{\gamma}}
\newcommand{\ugam}{\underline{\gamma}}
\newcommand{\olam}{\overline{\lambda}}
\newcommand{\ulam}{\underline{\lambda}}

\usepackage{dsfont,amssymb,amsmath, graphicx} 
\usepackage{amsfonts,dsfont,mathtools, mathrsfs} 
\usepackage{algorithm,algorithmicx,fancyhdr}
\usepackage[noend]{algpseudocode}

\usepackage{calrsfs}
\DeclareMathAlphabet{\pazocal}{OMS}{zplm}{m}{n}
\let\mathcal\undefined
\DeclareMathAlphabet{\mathcal}{OMS}{cmsy}{m}{n}

\newtheorem{definition}{Definition}
\newtheorem{assumption}{Assumption}
\newtheorem{corollary}{Corollary}
\newtheorem{theorem}{Theorem}
\newtheorem{remark}{Remark}

\newcommand{\Pb}{\mathbb{P}}
\newcommand{\Eb}{\mathbb{E}}

\newcommand{\R}{\mathds{R}}
\newcommand{\N}{\mathds{N}}

\newcommand{\cl}[1]{{\mathcal #1}}

\newcommand{\frk}[1]{{\mathfrak #1}}
\newcommand{\bs}[1]{{\boldsymbol #1}}
\newcommand{\tx}[1]{{\text #1}}

\newcommand{\col}{\text{col}}

\newcommand*{\QEDB}{\hfill\ensuremath{\square}}%

\usepackage{color}
\usepackage[dvipsnames]{xcolor}

\usepackage[
bookmarks=true,         
unicode=false,          
pdftoolbar=true,        
pdfmenubar=true,        
pdffitwindow=false,     
pdfstartview={FitV},    
pdftitle={},    
pdfauthor={},     
pdfsubject={},   
pdfcreator={},   
pdfproducer={}, 
pdfkeywords={}, 
pdfnewwindow=true,      
raiselinks	= true,
colorlinks=true,       
linkcolor=MidnightBlue,          
citecolor=Mahogany,        
filecolor=ForestGreen,      
urlcolor=ForestGreen,           
plainpages=false,
pdfpagelabels
]{hyperref}

\date{\today}

\usepackage{eurosym}
\usepackage{amstext} 
\DeclareRobustCommand{\officialeuro}{%
	\ifmmode\expandafter\text\fi
	{\fontencoding{U}\fontfamily{eurosym}\selectfont e}}

\usepackage[T1]{fontenc}     
\usepackage[utf8]{inputenc}	
\usepackage{graphicx, keyval, trig, url} 

\usepackage{epstopdf}

\usepackage{textcomp}
\usepackage{gensymb}
\usepackage[capitalise]{cleveref}
\usepackage{bm}   

\usepackage{blindtext}

\usepackage{cite} 

\usepackage{relsize}

\makeatletter

\def\subsection{\@startsection{subsection}{2}%
	\z@{.5\linespacing\@plus.5\linespacing}{.25\linespacing}%
	{\normalfont\bfseries}}

\def\subsubsection{\@startsection{subsubsection}{3}%
	\z@{.5\linespacing\@plus.5\linespacing}{.25\linespacing}%
	{\normalfont\itshape}}

\makeatletter
\def\@settitle{\begin{center}%
  \baselineskip14\p@\relax
	\normalfont\LARGE\scshape  
  \@title
  \end{center}%
}
\makeatother

\newcommand{\key}[1]{\textbf{Keywords.} #1}

\newcommand{\widthscale}{0.65\textwidth}

\title[Probabilistic Energy Management for BCC in STGs with Seasonal Storage Systems]
{Probabilistic Energy Management for Building Climate Comfort in Smart Thermal Grids with Seasonal Storage Systems} 

\author[V.\ Rostampour]{Vahab Rostampour}
\author[T.\ Keviczky]{Tam\'{a}s Keviczky}

\thanks{This research was supported by the Uncertainty Reduction in Smart Energy Systems (URSES) research program funded by the Dutch organization for scientific research (NWO) and Shell under the project Aquifer Thermal Energy Storage Smart Grids (ATES-SG) with grant number 408-13-030.\\
The authors are with Delft Center for Systems and Control, Delft University of Technology, Mekelweg 2, 2628 CD, Delft, The Netherlands. {\tt\small \{v.rostampour, t.keviczky\}@tudelft.nl}}%

\date{\today}

\begin{document} 

\begin{abstract}
This paper presents an energy management framework for building climate comfort (BCC) systems interconnected in a grid via aquifer thermal energy storage (ATES) systems in the presence of two types of uncertainty (private and common). 
ATES can be used either as a heat source (hot well) or sink (cold well) depending on the season. 
We consider the uncertain thermal energy demand of individual buildings as a private uncertainty source and the uncertain common resource pool (ATES) between neighbors as a common uncertainty source. 
We develop a large-scale stochastic hybrid dynamical model to predict the thermal energy imbalance in a network of interconnected BCC systems together with mutual interactions between their local ATES. 
We formulate a finite-horizon mixed-integer quadratic optimization problem with multiple chance constraints at each sampling time, which is in general a non-convex problem and difficult to solve. 
We then provide a computationally tractable framework by extending the so-called robust randomized approach and offering a less conservative solution for a problem with multiple chance constraints. 
A simulation study is provided to compare completely decoupled, centralized and move-blocking centralized solutions. 
We also present a numerical study using a geohydrological simulation environment (MODFLOW) to illustrate the advantages of our proposed framework.
\end{abstract}

\maketitle

\key{
	Smart Thermal Grids, 
	Building Climate Comfort Systems,
	Seasonal Storage Systems,
	ATES,
	Multiple Chance Constraints, 
	Probabilistic Robustness,
	Robust Randomized MPC.
}

\section{Introduction}
Global energy consumption has significantly increased due to the combined factors of increasing population and economic growth over the past few decades.
This increasing consumption highlights the necessity of employing innovative energy saving technologies.
Smart Thermal Grids (STGs) can play an important role in the future of the energy sector by ensuring a heating and cooling supply that is more reliable and affordable for thermal energy networks connecting various households, greenhouses and other buildings, which we refer to as agents.  
STGs allow for the adaptation to changing circumstances, such as daily, weekly or seasonal variations in supply and demand by facilitating each agent with smart thermal storage technologies.

Aquifer thermal energy storage (ATES) is a less well-known sustainable seasonal storage system that can be used to store large quantities of thermal energy in aquifers.  
Aquifers are underground porous formations containing water that are suitable for seasonal thermal energy storage.
It is especially suitable for climate comfort systems of large buildings such as offices, hospitals, universities, musea and greenhouses, see \cite{jaxa2016assessing}.
Most buildings in moderate climates have a heat shortage in winter and a heat surplus in summer. Where aquifers exist, this temporal discrepancy can be overcome by seasonally storing and extracting thermal energy into and out of the subsurface, enabling the reduction of energy usage and CO$_2$ emissions of climate comfort systems in buildings.
Fig.~\ref{fig_ates} depicts the operating modes of an ATES system for a single building.

\subsection*{Related Works}

There are various studies in literature related to buildings integrated into a smart grid \cite{razmara2016bilevel,taha2017buildings}. 
Modeling a building heating system  connected to a heat pump can be found in \cite{Michailidis2015}, an experimental model with  a  focus  on  heating,  ventilation,  and  air conditioning (HVAC) systems in \cite{stamatescu2016building}, using multi-HVAC systems in \cite{Satyavada2016}.
Models for building system dynamics together with HVAC controls are  typically linear \cite{wang2006parameter} for obvious computational  purposes.
For instance resistance and  capacitance circuit models, that represent heat transfer and thermodynamical  properties  of the building, are commonly used for building control studies \cite{dong2010integrated,maasoumy2014modeling,yahiaoui2006model}.
PID controllers for HVAC systems are widely used in many commercial buildings \cite{afram2014theory}.
Model  predictive  control  (MPC), on the other hand, has  received a lot of attention \cite{morocsan2011distributed,patel2016distributed,mirakhorli2016occupancy}, since it can handle large-scale dynamical systems subject to hard constraints, e.g., equipment limitations.
Using demand response for smart buildings \cite{qidemand}, MPC can be used in building climate comfort (BCC) problems \cite{ma2012model,oldewurtel2012use}.
MPC can overcome BCC problems even in decentralized or distributed setting and it is shown that has several advantages compared to PID controllers \cite{morocsan2011distributed,patel2016distributed,standardi2015economic}.

STGs have been studied implicitly in the context of micro combined heat and power systems, see \cite{ummenhofer2017improve}, or general smart grids, e.g., see  \cite{larsen2013distributed} and \cite{larsen2014distributed}.
Building heat demand with a dynamical storage tank was considered in \cite{van2013flexibility}, whereas in \cite{powell2013adaptive} an adaptive-grid model for dynamic simulation of thermocline thermal energy storage systems was developed. 
A deterministic view on STGs was studied by a few researchers \cite{rivarolo2013thermo}, \cite{lund20144th}, \cite{sameti2017optimization}.
STGs with uncertain thermal energy demands have been considered in \cite{farahani2016robust}, where a MPC strategy was employed with a heuristic Monte Carlo sampling approach to make the solution robust.	
A dynamical model of thermal energy imbalance in STGs with a probabilistic view on uncertain thermal energy demands was established in \cite{rostampour2016robust}, where a stochastic MPC with a theoretical guarantee on the feasibility of the obtained solution was developed.

\subsection*{Contributions}

ATES as a seasonal storage system has not, to the best of our knowledge, been considered in STGs.
In \cite{rostampour2016control} and \cite{rostampour2016building}, a dynamical model for an ATES system integrated in a BCC system has been developed. 
Following these studies, the first results toward developing an optimal operational framework to control ATES systems in STGs is presented here.  
In this framework, uncertain thermal energy demands are considered along with the possible mutual interactions between ATES systems, which may cause limited performance and reduced energy savings.
The main contributions of this paper are threefold:

a) We develop a novel large-scale stochastic hybrid dynamical model to predict the dynamics of thermal energy imbalance in STGs consisting of BCC systems with hourly-based operation and ATES as a seasonal energy storage system.
Based on our previous work in \cite{rostampour2016control} and \cite{rostampour2016building}, we extend an ATES system model to predict the amount of stored water and thermal energy. We first incorporate the ATES model into a BCC problem and then, formulate a large-scale STGs problem by taking into consideration the geographical coupling constraints between ATES systems.   
Using an MPC paradigm, we formulate a finite-horizon mixed-integer quadratic optimization problem with multiple chance constraints at each sampling time leading to a non-convex problem, which is difficult to solve.

b) We next propose a move-blocking control scheme to enable our stochastic MPC framework to handle long prediction horizons and  an hourly-based operation of the BCC systems together with a seasonal variation of desired optimal operation of the ATES system in a unified framework.
In practice, the BCC systems have an hourly-based operation and typically day-ahead planning compared to the ATES system that is based on a seasonal operation.
Using a fixed prediction horizon length, e.g., least common multiple of these two systems, may turn out to be computationally prohibitive, however also necessary in order to represent ATES interaction dynamics. 
The time scale discrepancy between the ATES system dynamics and BCC systems are explicitly accounted for in the developed MPC-based optimization formulation. 
Our proposed control strategy offers a long enough prediction horizon to prevent mutual interactions between ATES systems with much less computational time compared to a fixed  prediction horizon that is sampled densely (i.e., every hour).

c) We develop a computationally tractable framework to approximate a solution of our proposed MPC formulation based on our previous work in \cite{rostampour2016robust}.
In particular, we extend the framework in \cite{rostampour2016robust} to cope with multiple chance constraints which provides a more flexible approximation technique compared to the so-called robust randomized approach \cite{margellos2014road,margellos2013stochastic}, which is only suitable for a single chance constraint.
Our framework is closely related to, albeit different from, the approach of \cite{schildbach2013randomized}.
In \cite{schildbach2013randomized}, the problem formulation is convex and consists of an objective function with multiple chance constraints, in which the terms in objective and constraints are univariate.
In contrast, our problem formulation is mixed-integer and the objective function consists of separable additive components.

It is important to highlight that two major difficulties arising in stochastic hybrid MPC, namely recursive feasibility and stability, are not in the scope of this paper, and they are subject of our ongoing research work.
Thus, instead of analyzing the closed-loop asymptotic behavior, in this paper we focus on individual stochastic hybrid MPC problem instances from the optimization point of view and derive probabilistic guarantees for multiple chance constraints fulfillment.
A conference version of this paper was published in \cite{rostampour2017energymanagement}.
Additional contributions of the current manuscript: 
(i) we present a complete modeling of integrated BCC systems using ATES in STGs together with detailed MPC-based optimization formulation for each agent.
(ii) we provide a new control strategy, the so-called move-blocking control scheme in \Cref{sec_formulation}, to cope with the time-scale discrepancy between the ATES system dynamics and BCC systems.  
(iii) a complete proof of \Cref{cor_common_cons} and \Cref{thm_common} are presented in contrast to \cite{rostampour2017energymanagement};
(iv) results of an extensive simulation study is presented to compare all the proposed control schemes together with the move-blocking MPC scheme.

\begin{figure}
	\centering
	\includegraphics[width=\widthscale]{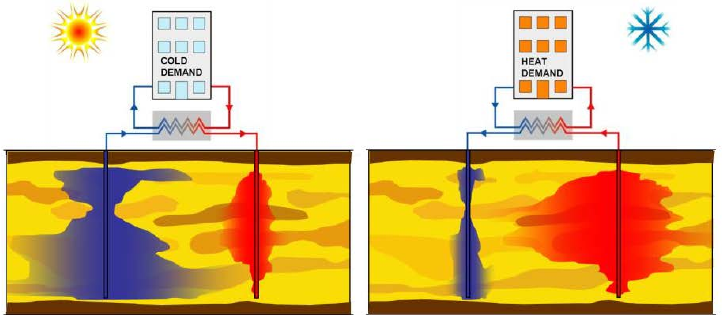}
	\caption{Operational modes of an ATES system during warm (left) and cold (right) seasons. Figure is taken from \cite{rostampour2016building}.}
	\label{fig_ates}
\end{figure}

\subsection*{Structure}

The layout of this paper is as follows.
Section \ref{sec_model} describes dynamics of an ATES system and a BCC system.
In Section \ref{sec_formulation}, we first formulate an energy management problem in a single agent, and then, extend it to a network of multiple agents.
We present three different setups, namely: one with completely decoupled agents, a centralized problem, and a move-blocking centralized problem formulation.
In Section \ref{sec_framework}, we develop a computationally tractable framework to solve these problems, and Section \ref{sec_sim} provides a simulation study with a comparison between these three different settings.  In addition, the numerical results obtained via a geohydrological simulation environment (MODFLOW) are shown.
Section \ref{final} concludes this paper with some remarks and future work.

\section*{Notation}

The following international system of units is used throughout the paper: 
Kelvin $[\tx{K}]$ and Celsius $[\cel]$ are the units of temperature,  Meter $[\tx{m}]$ is the unit of length, Hour $[\tx{h}]$ is the unit of time, Kilogram $[\text{kg}]$ is the unit of mass, Watt $[\tx{W}]$ is the unit of power, Joule $[\tx{J}]$, kiloWatt-hour $[\text{kWh}]$, and  MegaWatt-hour $[\text{MWh}]$ are the units of energy.

$\R, \R_{+}$ denote the real and positive real numbers, and $\N, \N_{+}$ the natural and positive natural numbers, respectively.
We operate within $n$-dimensional space $\R^n$ composed by column vectors $u, v \in \R^n$.
The Cartesian product over $n$ sets $\cl{X}_1, \cdots, \cl{X}_n$ is given by: $\prod_{i=1}^n \cl{X}_i = \cl{X}_1 \times \cdots \times \cl{X}_n = \{(x_1, \cdots, x_n) \, : \, x_i \in \cl{X}_i\}$.
The cardinality of a set $\cl{A}$ is shown by $|\cl{A}| = A$.

Given a metric space $\Delta$, and $\Pb$ a probability measure defined over $\Delta$, its Borel $\sigma$-algebra is denoted by $\frk{B}(\Delta)$.
Throughout the paper, measurability always refers to Borel measurability.
In a probability space $(\Delta, \frk{B}(\Delta),\Pb)$, we denote the $N$-Cartesian product set of $\Delta$ by $\Delta^N$ with the respective product measure by $\Pb^N$.

\section{System Dynamics Modeling}\label{sec_model}

In this section, we first develop a mathematical model for an ATES system dynamics as a single seasonal storage system. 
We then describe the steady-state dynamical model for a building to capture its thermal energy demand profile during heating and cooling modes based on our previous work in \cite{rostampour2017hpc}.
We finally present the BCC system where we introduce the so-called thermal energy imbalance dynamics.
Using the thermal energy imbalance dynamics of a BCC system, we integrate an ATES system into the building thermal energy production unit which consists of a boiler, a chiller, a heat pump, a heat exchanger and storage tanks as the heating and cooling modes equipment.

\subsection{Seasonal Storage Systems}

Consider an ATES system consisting of warm and cold wells to store warm water during warm season and cold water during cold season, respectively.
Each well can be described as a single thermal energy storage where the amount of stored energy is proportional to the temperature difference between stored water and aquifer ambient water. 
Stored thermal energy from the last season is going to be used for the current season and so forth. 
Depending on the season, the operating mode (heating or cooling) of an ATES system changes, by reversing the direction of water flow between wells
as it is shown in Fig.~\ref{fig_ates}.

During a cold season, for heating purposes, the direction of water is from the warm well to the cold well through a heat exchanger to extract the stored thermal energy from the water. The return water is cooled down and stored in the cold well.  
This procedure is opposite during a warm season for cooling purposes of the BCC system.
An ATES system can be characterized by some physically meaningful parameters.
The most relevant features that can describe the status of an ATES system for the purpose of optimal control is the stored volume of water together with the thermal energy content in each well. 
A free manipulated variable in this setting is the pump flow rate that is used to circulate water from one well to the other through a heat exchanger.

We therefore define the states that can describe the ATES system dynamics to be the volume of water, $\tx{V}_{\tx{a},k}^{h} \, [\tx{m}^3]$, $\tx{V}_{\tx{a},k}^{c} \, [\tx{m}^3]$, and the thermal energy content, $\tx{S}_{\tx{a},k}^{h} \, [\text{Wh}], \tx{S}_{\tx{a},k}^{c} \, [\text{Wh}]$, of warm and cold wells. 
The superscripts "$h$" and "$c$" refer to the heating and cooling operating modes of an ATES system, respectively, and the subscript "a" denotes the ATES system variables.
Consider the following first-order difference equations as ATES system model dynamics:
\begin{subequations}\label{ates_dyn}
	\begin{align}
	\tx{V}_{\tx{a},k+1}^h &= \tx{V}_{\tx{a},k}^h -  \tau (u_{\tx{a},k}^h - u_{\tx{a},k}^c) \, , \label{ates_dynVh} \\
	\tx{V}_{\tx{a},k+1}^c &= \tx{V}_{\tx{a},k}^c +  \tau (u_{\tx{a},k}^h - u_{\tx{a},k}^c) \, , \label{ates_dynVc} \\
	\tx{S}_{\tx{a},k+1}^h &= \eta_\tx{a,k} \, \tx{S}_{\tx{a},k}^h -  \tau (h_{\tx{a},k}^h - h_{\tx{a},k}^c) \, , \label{ates_dynSh} \\
	\tx{S}_{\tx{a},k+1}^c &= \eta_\tx{a,k} \, \tx{S}_{\tx{a},k}^c +  \tau (c_{\tx{a},k}^h - c_{\tx{a},k}^c) \, , \label{ates_dynSc}
	\end{align}
\end{subequations}
where $\eta_{\tx{a,k}} \in (0,1)$ is a lumped coefficient of thermal energy losses in aquifers,
$u_{\tx{a},k}^h \, [\tx{m}^3 \tx{h}^{-1}]$, and $u_{\tx{a},k}^c \, [\tx{m}^3 \tx{h}^{-1}]$ are control variables corresponding to the pump flow rate of ATES system during heating and cooling modes at each sampling time $k = 1, 2, \cdots$, respectively, with $\tau\, [\text{h}]$ as the sampling period.
$u_{\tx{a},k}^h$ circulates water from warm well to cold well, whereas $u_{\tx{a},k}^c$ takes water from cold well and injects into warm well of ATES system, during heating modes and cooling modes of the BCC system, respectively. 
The variables $h_{\tx{a},k}^h \, [\tx{W}], c_{\tx{a},k}^h \, [\tx{W}]$ denote the thermal power that is extracted from warm well and injected into cold well of ATES system during heating mode of BCC system, respectively. 
The variables $c_{\tx{a},k}^c \, [\tx{W}], h_{\tx{a},k}^c \, [\tx{W}]$ are the thermal power that is extracted from cold well and injected into warm well of ATES system during cooling mode of BCC system, respectively.
These variable are defined by:
\begin{align}\label{ATESvaryingVars}
\begin{cases}
h_{\tx{a},k}^h = \alpha_{h,k} \, u_{\tx{a},k}^h \\
c_{\tx{a},k}^h = \alpha_{c,k} \, u_{\tx{a},k}^h 
\end{cases},
\begin{cases}
c_{\tx{a},k}^c = \alpha_{c,k} \, u_{\tx{a},k}^c \\ 
h_{\tx{a},k}^c = \alpha_{h,k} \, u_{\tx{a},k}^c  
\end{cases},
\end{align}
where $\alpha_{h,k} = \rho_w \, c_{pw} \,(\tx{T}_{\tx{a},k}^{h} - \tx{T}_{\tx{a},k}^{\text{amb}})$, and $\alpha_{c,k} = \rho_w \, c_{pw} \, (\tx{T}_{\tx{a},k}^{\text{amb}} - \tx{T}_{\tx{a},k}^{c})$ are the thermal power coefficients of warm and cold wells, respectively.
The parameters $\rho_w$ $[\text{kg} \tx{m}^{-3}]$, $c_{pw}$ $[\tx{J} \text{kg}^{-1} \tx{K}^{-1}]$ are density and specific heat capacity of water, respectively.
$\tx{T}_{\tx{a},k}^{h}$ $[\tx{K}]$, $\tx{T}_{\tx{a},k}^{c}$ $[\tx{K}]$, and $\tx{T}_{\text{aq},k}^{\text{amb}}$ $[\tx{K}]$ denote the temperature of water inside warm well, cold well and the ambient aquifer, respectively.
We also define $h_{\tx{a},k} \, [\text{Wh}]$, and $c_{\tx{a},k} \, [\text{Wh}]$ to be the amount of thermal energy that can be delivered to the building during heating and cooling modes, respectively, as follows:
\begin{align}\label{energy_delivered_ates}
\begin{cases}
h_{\tx{a},k} = \alpha_k \,\tau\, u_{\tx{a},k}^h \\
c_{\tx{a},k} = \alpha_k \,\tau\, u_{\tx{a},k}^c 
\end{cases}, 
\end{align}
where $\alpha_k = \alpha_{h,k} + \alpha_{c,k}$ is the total thermal power coefficient.
In previous work \cite{rostampour2016control}, we have also developed a control-oriented model for the integrated ATES into BCC system where we consider the dynamical behavior of the volume and temperature of water in each well of ATES system. 

\begin{figure}
	\centering
	\includegraphics[width=\widthscale]{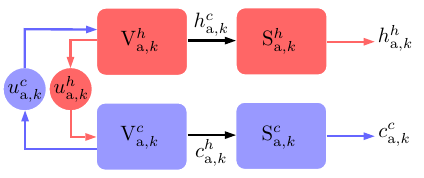}
	\caption{
		Operational block diagram of ATES system showing the relation between an ATES system variables. Operating modes during cold and warm seasons are shown via red color and blue color, respectively.
	}
	\label{fig_ates_block}
\end{figure}

Let us now discuss the dynamics of ATES system in \eqref{ates_dyn}.
\Cref{ates_dynVh,ates_dynVc,ates_dynSh,ates_dynSc} describe the evolution of water volume and the thermal energy content in warm and cold wells, respectively. 
During cold seasons, for the heating purpose of BCC system, the amount of $\tau u_{\tx{a},k}^h$ volume of warm water from the warm well is extracted to provide $\tau h_{\tx{a},k}$ amount thermal energy, and meanwhile, the amount of $\tau c_{\tx{a},k}^h$ thermal energy is stored in the cold well of ATES system.
As for the cooling purpose of BCC system during warm seasons, the amount of $\tau u_{\tx{a},k}^c$ volume of cold water from the cold well is extracted to provide $\tau c_{\tx{a},k}$ amount thermal energy, while the amount of $\tau h_{\tx{a},k}^c$ thermal energy is injected in the warm well of ATES system.
The operations of the ATES system is visualized in Fig.~\ref{fig_ates_block}, which represents the relation between the variables of ATES system. 
Operating modes during cold and warm seasons are shown via red color and blue color, respectively.
The following assumption is made due to the existing operational practice, and it is not restrictive for our proposed model.
\begin{assumption}
	There is either no operation or only one operating mode active in ATES systems, which leads to either both control variables being zero or only one control variable being nonzero at any time instant:
	\begin{align*}
	u_{\tx{a},k}^h \, u_{\tx{a},k}^c = 0 \,, \quad k = 1, 2, \cdots \ .
	\end{align*}
\end{assumption}

The dynamics of ATES system in \eqref{ates_dyn} can be also written in a more compact format for each agent $i \in \{1, \cdots, N\}$:
\begin{align}\label{atesi}
x_{i,k+1}^\tx{a} = a_{i,k}^\tx{a} x_{i,k}^\tx{a} + b_{i,k}^\tx{a} u_{i,k}^\tx{a} \ , 
\end{align}
where 
$x_{i,k}^\tx{a} = \begin{bmatrix} \tx{V}_{\tx{a},k}^h & \tx{V}_{\tx{a},k}^c & \tx{S}_{\tx{a},k}^h & \tx{S}_{\tx{a},k}^c \end{bmatrix}^\top \in\R^4$ denotes the state vector, 
$u_{i,k}^\tx{a} = \begin{bmatrix} u_{\tx{a},k}^h & u_{\tx{a},k}^c \end{bmatrix}^\top \in\R^2$ is the control vector, and $a_{i,k}^\tx{a}, b_{i,k}^\tx{a}$ can be obtained via \eqref{ates_dyn}.
Note that there are some operational constraints on the ATES control variable as well,
\begin{subequations}\label{atesUcons}
	\begin{align}
	& u_{\tx{a}}^{\min} \leq u_{\tx{a},k}^{h} \leq u_{\tx{a}}^{\max}\,, \\
	& u_{\tx{a}}^{\min} \leq u_{\tx{a},k}^{c} \leq u_{\tx{a}}^{\max}\,, 
	\end{align}
\end{subequations}
where $u_{\tx{a}}^{\min}, u_{\tx{a}}^{\max}$ represent the minimum and maximum pump flow rate of ATES system, respectively.

The proposed model for an ATES system in \eqref{atesi} is a linear time-varying discrete-time system, due to the variation of the temperatures in both wells and the ambient aquifer \eqref{ATESvaryingVars}.
In \Cref{model_bcc}, we will integrate \eqref{atesi} into a BCC system dynamics.

\subsection{Thermal Energy Demand Profile}\label{thermal_energy}

A dynamical model of building thermal energy demand was developed in our previous work \cite{rostampour2017hpc} to determine the thermal energy demand of a building at each sampling time $k$, considering the desired indoor air temperature and the outside weather conditions.  
We refer to the BCC system that determines the level of thermal energy demand $\tx{Q}_{\tx{d},k}^\tx{B}$ $[\text{Wh}]$ at each sampling time $k$ via 
\begin{align}\label{QB}
\tx{Q}_{\tx{d},k}^\tx{B} = f_{\tx{B}}(p^\tx{B}_s, \tx{T}_{\text{des},k}^\tx{B}, \vartheta_{k}) \,,
\end{align}
where $p^\tx{B}_s$ corresponds to a parameter vector of the building characteristics, $\tx{T}_{\text{des},k}^\tx{B}$ $[\cel]$ is the desired indoor air temperature of the building, and 
$\vartheta_{k} = [\tx{T}_{o,k}^\tx{B}, \tx{I}_{o,k}, \tx{v}_{o,k}, \tx{Q}_{p,k}, \tx{Q}_{e,k}] \in \R^{5}$ is a vector of uncertain variables  that contains the outside air temperature, the solar radiation, the wind velocity, the thermal energy produced due to occupancy by people, and electrical devices, and lighting inside the building, respectively.
This yields the building thermal energy demand that takes into account the overall building effects, e.g., zones, walls, humans and non-human thermal energy sources with the outside uncertain weather conditions.
Since we are mainly interested in capturing the variation of thermal energy demand w.r.t. the outside air temperature $\tx{T}_{o,k}^\tx{B}$, the uncertain variable $\vartheta_{k}$, is assigned to $\tx{T}_{o,k}^\tx{B}$, and the rest of the variables are fixed to their nominal (forecast) values at each sampling time $k$.

The operating modes (heating or cooling) of BCC system are determined based on the sign of $\tx{Q}_{\tx{d},k}^\tx{B}$ at each sampling time $k$. 
The variable $\tx{Q}_{\tx{d},k}^\tx{B}$ with positive and negative signs, represents the thermal energy demand during heating mode and the building surplus thermal energy during cooling mode, respectively. 
$\tx{Q}_{\tx{d},k}^\tx{B}$ is zero represents the comfort mode of building, and thus, in such a case no heating or cooling is requested.
We also distinguish between the thermal energy demand of building during heating mode $h_{\tx{d},k}$, and cooling mode $c_{\tx{d},k}$, using the relation: $\tx{Q}_{\tx{d},k}^\tx{B} = h_{\tx{d},k} - c_{\tx{d},k}\,.$

The following technical assumption is necessary for the measurability of the uncertainty. 

\begin{assumption}
	The mapping from the uncertain variable $\vartheta_{k}$ to the thermal energy demand $\tx{Q}_{\tx{d},k}^\tx{B}$ is a measurable function \eqref{QB}, so that $\tx{Q}_{\tx{d},k}^\tx{B}$ can be viewed as a random variable on the same probability space as $\vartheta_{k}$. 
	Moreover, the thermal energy demand at each sampling time $k$ 
	can be either zero (no thermal energy demand) or only for heating $h_{\tx{d},k}$ (cooling $c_{\tx{d},k}$) mode, similarly to the operating modes of the ATES system, which leads to: 
	\begin{align*}
	h_{\tx{d},k} \ c_{\tx{d},k} = 0 \ , \ k = 1, 2, \cdots \ .
	\end{align*}  
\end{assumption}

\begin{figure}
	\centering
	\includegraphics[width=\widthscale]{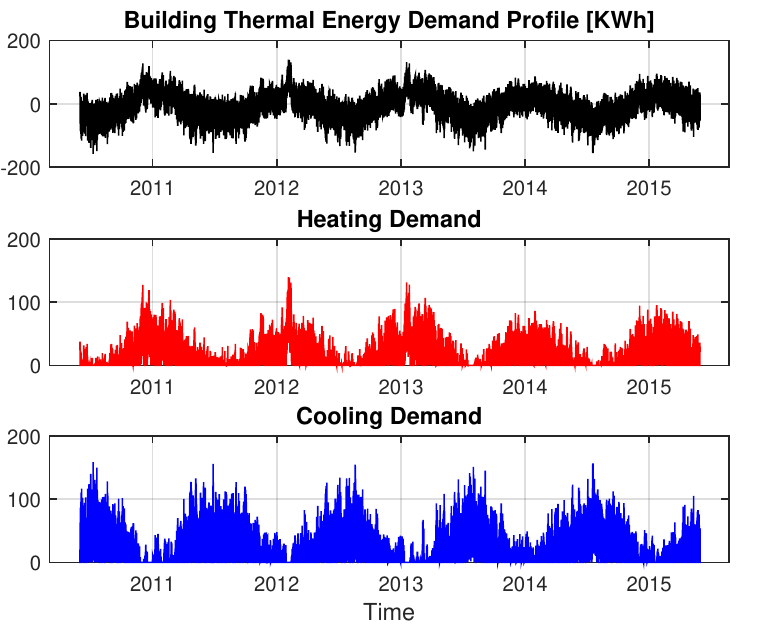}
	\caption{Thermal energy demand profile of a building during 2010-2015 with respect to the outside registered weather data in The Netherlands \cite{wicak2016distributed}. The black line shows $\tx{Q}_{\tx{d},k}^\tx{B}$, the red line is related to the thermal energy demand during heating mode $h_{\tx{d},k}$ and the blue line corresponds to the thermal energy surplus during cooling mode $c_{\tx{d},k}$.}
	\label{fig_demand_profile}
\end{figure}

Fig.~\ref{fig_demand_profile} shows the thermal energy demand profile of a building for the last five years with respect to the outside registered weather data in The Netherlands. 
The top panel in Fig.~\ref{fig_demand_profile} depicts $\tx{Q}_{\tx{d},k}^\tx{B}$ as the result of \eqref{QB}, whereas the middle and bottom panels show the thermal energy demand during heating mode $h_{d,k}$ and the thermal energy surplus during cooling mode $c_{d,k}$, respectively.

\subsection{Building Climate Comfort Systems}\label{model_bcc}

\begin{figure}
	\centering
	\includegraphics[width=\widthscale]{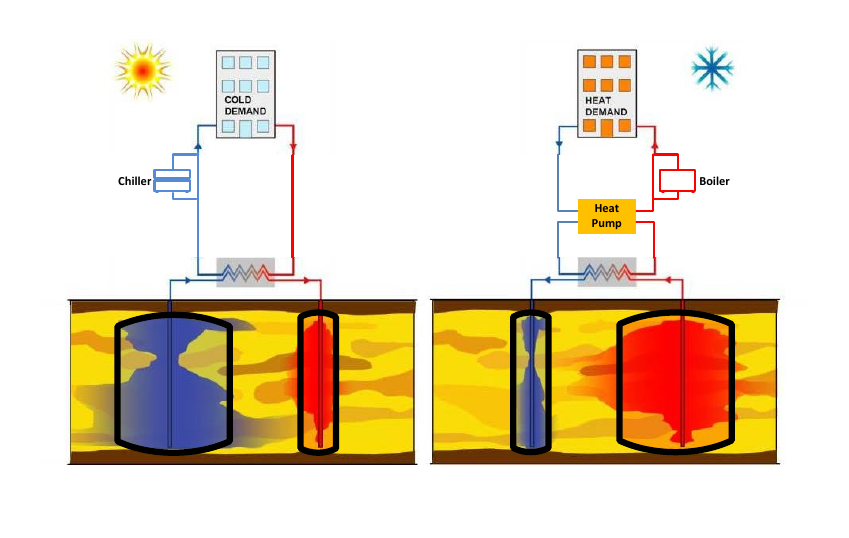}
	\caption{Heating and cooling operating modes of BCC system with an ATES system during warm (left) and cold (right) seasons.} 
	\label{fig_hcnets}
\end{figure}

Consider a single agent (i.e., building) $i \in \{1, \cdots, N\}$ that is facilitated with a boiler, a heat pump, a storage tank for the heating mode, and a chiller, a storage tank for the cooling mode together with an ATES system that is available for both operating modes (see Fig.~\ref{fig_hcnets}). We now focus on the modeling of energy balance for the BCC system. 

Define two vectors of control variables during heating and cooling modes in each agent $i$ at each sampling time $k$, to be
\begin{align*}
u_{i,k}^{h} = \begin{bmatrix} h_{\boi,k}^{} & h_{\im,k} \end{bmatrix}^\top \in \R^{2} \,,\,
u_{i,k}^{c} = \begin{bmatrix} c_{\chiv,k}^{} & c_{\im,k} \end{bmatrix}^\top \in \R^{2} \,.
\end{align*}
The variables $h_{\boi,k}$, $c_{\chiv,k}$, $h_{\im,k}\,$, and $c_{\im,k}$ denote the production of boiler, chiller, the imported energies from external parties during heating and cooling modes, respectively.
We consider boiler and chiller operating limits that constrain their production within a certain bound for cost effective maintenance of such equipment. 
Define $v_{\boi,k}\in \{0,1\}$ and $v_{\chiv,k} \in \{0,1\}$ to be two binary variables to decide about the ON/OFF status boiler and chiller, respectively.
Consider now to the following conditional situations:
\begin{subequations}
	\label{status_equip}
	\begin{align}
	\text{boiler}\,&:\,
	\begin{cases}
	v_{\boi,k} = 1  \quad & h_{\boi}^{\min} \leq h_{\boi,k}^{} \leq h_{\boi}^{\max}\\
	v_{\boi,k} = 0   \quad & \quad\ \, \text{otherwise}
	\end{cases}, \\
	\text{chiller}\,&:\,
	\begin{cases}
	v_{\chiv,k} = 1 \quad & c_{\chiv}^{\min} \leq c_{\chiv,k}^{} \leq c_{\chiv}^{\max} \\
	v_{\chiv,k} = 0 \quad & \quad\ \, \text{otherwise} 
	\end{cases}\,,
	\end{align}
\end{subequations}
where $h_{\boi}^{\min}$, $h_{\boi}^{\max}$, $c_{\chiv}^{\min}$, $c_{\chiv}^{\max}$ denote the minimum and maximum capacity of thermal energy production of boiler and chiller, respectively.

We define two variables to capture the thermal energy imbalance errors during heating mode $x_{i,k}^h \in \R$, and an imbalance error of the cooling mode $x_{i,k}^c \in \R$.
They are related to the difference between the level of the storage tank with the forecasted thermal energy demand, $h_{\tx{d},k}^{\tx{f}}$, $c_{\tx{d},k}^{\tx{f}}$, during heating and cooling modes, respectively, which are formally defined using the following relations:
\begin{subequations}
	\label{xdef} 
	\begin{align}
	x_{i,k}^h = h_{\tx{s},k} - h_{\tx{d},k}^{\tx{f}} \ , \\
	x_{i,k}^c = c_{\tx{s},k} - c_{\tx{d},k}^{\tx{f}} \ . 
	\end{align}
\end{subequations}
Herein, $h_{\tx{s},k}$, and $c_{\tx{s},k}$ represent the level of storage tank during heating and cooling modes, respectively, and obey the following dynamics:
\begin{align}\label{bufferdyn}
h_{\tx{s},k+1} &= \eta_{\tx{s},k}^h x_{i,k}^h + \eta_{\tx{s},k}^h \left(h_{\boi,k} + h_{\im,k} + \alpha_{\text{hp},k} h_{\tx{a},k} \right) , \\
c_{\tx{s},k+1} &= \eta_{\tx{s},k}^c x_{i,k}^c + \eta_{\tx{s},k}^c \left(c_{\chiv,k} + c_{\im,k} + c_{\tx{a},k} \right) , 
\end{align}
where $\alpha_{\text{hp},k} = \text{COP}_k (\text{COP}_k-1)^{-1}$ is related to the effect of the heat pump during heating mode and COP$_k$ stands for the coefficient of performance of heat pump at each sampling time $k$.
The parameters $\eta_{\tx{s},k}^h \,, \eta_{\tx{s},k}^c \in (0,1)$ denote the thermal loss coefficients due to inefficiency of storage tank during heating and cooling modes, respectively.
The variables $h_{\tx{a},k}$ and $c_{\tx{a},k}$ are defined in \eqref{energy_delivered_ates} and are related to the ATES system model.
It is important to note that $h_{\tx{a},k}$ and $c_{\tx{a},k}$ are dependent on the pump flow rates $u_{\tx{a},k}^h$ and $u_{\tx{a},k}^c$ of the ATES system during heating and cooling modes of the BCC system, respectively.
We now substitute $h_{\tx{s},k}$, and $c_{\tx{s},k}$ as in \eqref{bufferdyn} into \eqref{xdef} to derive the dynamical behavior of the thermal energy imbalance $x_{i,k}^h$ and $x_{i,k}^c$ that are given by
\begin{subequations}
	\label{imbhc}
	\begin{align}
	x_{i,k+1}^h &= a_{i,k}^h x_{i,k}^h + b_{i}^h u_{i,k}^h + b_{i,k}^{\tx{a},h} u_{i,k}^\tx{a} + c_{i,k}^h w_{i,k}^h \ , \\
	x_{i,k+1}^c &= a_{i,k}^c x_{i,k}^c + b_{i}^c u_{i,k}^c + b_{i,k}^{\tx{a},c} u_{i,k}^\tx{a} + c_{i,k}^c w_{i,k}^c \ , 
	\end{align}
\end{subequations}
where 
$a_{i,k}^h = \eta_{\tx{s},k}^h$\,, 
$a_{i,k}^c = \eta_{\tx{s},k}^c$\,, 
$b_{i,k}^h = \begin{bmatrix} \eta_{\tx{s},k}^h &\eta_{\tx{s},k}^h \end{bmatrix}$\,, 
$b_{i,k}^c =  \begin{bmatrix} \eta_{\tx{s},k}^c & \eta_{\tx{s},k}^c \end{bmatrix}$\,, 
$b_{i,k}^{\tx{a},h}  = \begin{bmatrix}\eta_{\tx{s},k}^h \alpha_{\text{hp},k} \, \alpha_k & 0 \end{bmatrix}$\,, 
$b_{i,k}^{\tx{a},c}  = \begin{bmatrix}\eta_{\tx{s},k}^c \, \alpha_k & 0 \end{bmatrix}$\,, 
$c_{i}^h = -1\,,$ and $c_{i}^c = -1$.
The variables $w_{i,k}^h = h_{\tx{d},k+1}^{\tx{f}}$ and $w_{i,k}^c = c_{\tx{d},k+1}^{\tx{f}}$ refer to the forecast of thermal energy demand during heating and cooling modes in the next time step, respectively. 
The only uncertain variable in each agent $i$ at each sampling time $k$ is considered to be the deviation of actual thermal energy demand from its forecast value as defined in \Cref{thermal_energy}, and therefore, $w_{i,k}^h$ and $w_{i,k}^c$ represent uncertain parameters.

Consider now the system dynamics for each agent $i$ by concatenating the thermal energy imbalance errors during heating and cooling modes \eqref{imbhc} together with the state vector of the ATES system \eqref{atesi} as follows:
\begin{align}\label{imb}
x_{i,k+1} = a_{i,k} x_{i,k} + b_{i,k} u_{i,k} + c_{i,k} w_{i,k} \ , 
\end{align}
where 
$x_{i,k} = \begin{bmatrix} x_{i,k}^{h \top} & x_{i,k}^{c \top} & x_{i,k}^{\tx{a} \top} \end{bmatrix}^\top \in \R^{6}$ denotes the state vector, 
$u_{i,k} = \begin{bmatrix} u_{i,k}^{h \top} & u_{i,k}^{c \top} & u_{i,k}^{\tx{a}  \top} \end{bmatrix}^\top \in \R^{6}$ is the control vector,
and $w_{i,k} = \begin{bmatrix} w_{i,k}^h & w_{i,k}^c \end{bmatrix}^\top \in \mathcal{W}_{i,k}\subseteq\R^{2}$ is the uncertainty vector such that $\mathcal{W}_{i,k}$ is an unknown uncertainty set.
The system matrices $a_{i,k}\,, b_{i,k}\,, c_{i,k}$ can be readily derived from their definitions and we omit them in the interest of space.

The proposed model for a BCC system in \eqref{imb} is a stochastic hybrid linear time-varying discrete-time system. 
It is important to note that the hybrid nature of \eqref{imb} is due to the fact that each equipment (boiler and chiller) can be either ON or OFF as in \eqref{status_equip} depending on heating and cooling modes of the building.
This possibility therefore changes the proposed thermal energy imbalance error dynamics \eqref{imbhc}.

In order to provide a desired thermal comfort for each  BCC system in the following section, we will develop a control framework based on the MPC paradigm where \eqref{imb} is used to predict the thermal energy imbalance error dynamics together with the ATES system dynamics for each agent $i \in \{1, \cdots, N\}$, and then, extend this to a network of interconnected BCC systems.   
Moreover, we will provide a solution method to overcome an important challenge of the network of BCC systems due to the spatial distribution of ATES systems.  
An important remark is that the variations of system parameters in the proposed dynamical model \eqref{imb} evolve on a much slower time-scale compared to the system dynamics and, therefore, we consider the system dynamics  \eqref{imb} to be time-invariant in the following parts.
It is worth mentioning that our proposed control technique in this paper can be easily extended to cope with time-varying parameters by considering them as multiplicative uncertainty sources, see e.g., \cite{rostampour2017distributedsmpc}.

\section{Energy Management Problem}\label{sec_formulation}

In this section, we formulate an optimization problem for heating and cooling modes of the BCC system integrated with ATES which we refer to as a single agent energy management problem. 
We then extend the single agent problem to a network with multiple agents that can be producers and consumers of thermal energy in a STGs setting. 
In such a setting, there might be some unwanted overlap (mutual interactions) between the stored water in the wells of neighboring ATES systems (see Fig.~\ref{fig_atessg}) through aquifers.
Such an unwanted mutual interactions between warm and cold wells, clearly, reduce the energy efficiency of the ATES systems.
The goal of the agents is to match the local consumption and production and to avoid mutual interactions between their ATES systems in the network and thereby improve energy efficiency.

\subsection{Energy Balance in Single Agent System}

Consider an MPC problem with a finite prediction horizon $N_h$ for each agent $i\in\{1, \cdots, N\}$, and introduce the subscript $t$ in our notation to characterize the value of the planning quantities for a given time $t \in \mathcal{T}$, where the set of predicted time steps is denoted by $\mathcal{T}:=\{k, k+1, \cdots, k+N_h-1\}$.
Using the subscript $t|k$, we refer to the $t$ time step prediction of variables at the simulation time step $k$. 

Define $v_{i,t|k} = \begin{bmatrix} v_{\boi,t|k} & v_{\chiv,t|k} \end{bmatrix}^\top \in \{0,1\}^2$ as a vector of binary variables to decide about the ON/OFF status of boiler and chiller in each agent $i\in\{1, \cdots, N\}$.
We also take into account the startup cost of boiler and chiller using 
$c_{i,t|k}^{\text{su}} = \begin{bmatrix} c_{\boi,t|k}^{\text{su}} & c_{\chiv,t|k}^{\text{su}} \end{bmatrix}^\top$ and add $c_{i,t|k}^{\text{su}}$ into the control decision variables
$u_{i,t|k} = \begin{bmatrix} u_{i,t|k}^{h \top} & u_{i,t|k}^{c \top} & u_{i,t|k}^{\tx{a}  \top} & c_{i,t|k}^{\text{su} \top}\end{bmatrix}^\top \in \R^{8}$
for each agent $i$ at each time step $t|k$.

The goal of each agent $i$ is to map the local thermal energy supply of production units to the local thermal energy demand of BCC system. 
Our goal thus is to formulate an optimization problem to find the control input $u_{i,t|k}$ for each agent $i$ such that the thermal energy imbalance errors stay as a small as possible at minimal production cost and to satisfy physical constraints of heating and cooling modes equipment at each sampling time $k$.
We therefore associate a quadratic cost function with each agent $i$ at each prediction time step $k$ as follows: 
\begin{align}
J_i(x_{i,t|k}, u_{i,t|k}) = x_{i,t|k}^\top \, Q_i \, x_{i,t|k} + u_{i,t|k}^\top \, R_i \, u_{i,t|k} \ , 
\end{align}
where $Q_i = \text{diag}\left(\begin{bmatrix} q_i^h & q_i^c & \bs{0}_{1\times 4} \end{bmatrix}\right) \in \R^{6\times 6}$ is a weighting matrix coefficient of thermal energy imbalance errors, $R_i = \text{diag}\left(r_{i}\right)\in \R^{8\times 8}$ indicates a diagonal matrix with the cost vector $r_{i}$ on its diagonal, and
$r_{i}$ is defined as 
\begin{align*}
r_{i} = \begin{bmatrix} r_{\boi} & r_{\im}^h & r_{\chiv} & r_{\im}^c & r_{\tx{a}}^h & r_{\tx{a}}^c & 1 & 1 \end{bmatrix}^{\top} \in \R^{8} \ ,
\end{align*}
where 
$r_{\boi}\, (r_{\chiv})$ represents the cost of natural gas that is used by boiler (chiller), 
$r_{\im}^h\, (r_{\im}^c)$ denotes the cost of imported thermal energy from an external party during heating (cooling) mode, and 
$r_{\tx{a}}^h\, (r_{\tx{a}}^c)$ corresponds to the pumping electricity cost of ATES system to extract the required thermal energy during heating (cooling) modes.
The other entries of $r_{i}$ represent the start-up costs. 
The proposed cost function consists of two main parts which leads to the regulation of imbalance errors to zero at minimal production cost together with minimum energy balance error of ATES system in each agent $i$. 
The reason for introducing a cost function in this form is that from a computational perspective quadratic cost functions are motivated by convexity and differentiability arguments.
Note that the cost function $J_i(\cdot)$ is a random variable due to the uncertain state variables, and thus, we consider $\Eb_{}\left[J_i(\cdot)\right]$ to obtain a deterministic cost function.

We are now in a position to formulate a finite-horizon stochastic hybrid control problem as the local energy management problem for each agent $i\in\{1, \cdots, N\}$ using the following chance-constrained mixed-integer optimization problem:
\begin{subequations}
	\label{opt1}
	\begin{align}
	\underset{\substack{\{u_{i,t|k}, v_{i,t|k}\}_{t\in\mathcal{T}}}}{\text{min}} \quad &  \sum_{t\in\mathcal{T}}^{} \ \Eb\big[J_i(x_{i,t|k}, u_{i,t|k})\big]  \label{costopt1} \\
	\text{subject to} \qquad
	& \label{startup} c_{i,t|k}^{su} \geq \Lambda^{su} (v_{i,t|k} - v_{i,t-1|k}) \geq 0 \,, \forall t\in\mathcal{T}\\
	& \label{boilim} v_{\boi,t|k}  h_{\boi}^{\min} \leq h_{\boi,t|k} \leq h_{\boi}^{\max} v_{\boi,t|k}, \forall t\in\mathcal{T}\\
	& \label{chilim} v_{\chiv,t|k}  c_{\chiv}^{\min} \leq c_{\chiv,t|k} \leq c_{\chiv}^{\max} v_{\chiv,t|k}\,, \forall t\in\mathcal{T}\\
	& \label{him} h_{\im}^{\min} \leq h_{\im,t|k} \leq h_{\im}^{\max}\,, \forall t\in\mathcal{T}\\
	& \label{cim} c_{\im}^{\min} \leq c_{\im,t|k} \leq c_{\im}^{\max}\,, \forall t\in\mathcal{T}\\
	& \label{huaq} u_{\tx{a}}^{\min} \leq u_{\tx{a},t|k}^{h} \leq u_{\tx{a}}^{\max}\,, \forall t\in\mathcal{T}\\
	& \label{cuaq} u_{\tx{a}}^{\min} \leq u_{\tx{a},t|k}^{c} \leq u_{\tx{a}}^{\max}\,, \forall t\in\mathcal{T}\\
	& \label{probimb} \Pb\big\{ x_{i,t+1|k} \geq 0 \,\big\vert\, x_{i,t|k} \,,\,\forall t\in\mathcal{T} \big\} \geq 1-\varepsilon_i \,,\ \forall \{w_{i,t|k}\}_{t\in\mathcal{T}}\in\mathcal{W}_{i}  \,,
	\end{align}
\end{subequations}
where $\Lambda^{su}$ is a diagonal matrix including the startup costs of boiler and chiller on the diagonal, 
$h_{\im}^{\min}$, $h_{\im}^{\max}$, $c_{\im}^{\min}$, $c_{\im}^{\max}$ are the minimum and maximum capacity of thermal energy production for each external party during heating and cooling modes, respectively. 
$\varepsilon_i \in (0,1)$ is the admissible constraint violation parameter.
Note that $\mathcal{W}_{i}$ represents the Cartesian product of $\mathcal{W}_{i,t|k}$ for all $t\in\mathcal{T}$.

In order of appearance, the constraints have the following meaning.
Constraint \eqref{startup} captures the status change of boiler and chiller (from OFF to ON). 
Note that the status change from ON to OFF never appears in the cost function due to the positivity constraint of $c_{i,t|k}^{su} \geq 0$.     
\eqref{boilim}, \eqref{chilim}, \eqref{him}, \eqref{cim}, \eqref{huaq}, \eqref{cuaq} impose box constraints (capacity limitations) on their variables.
In the given lower and upper bounds of both constraints \eqref{boilim} and \eqref{chilim}, there are multiplications with binary variables which enforce the status change of boiler and chiller, respectively. 
Constraint \eqref{probimb} ensures probabilistically feasible trajectories of the thermal energy imbalance errors for in each agent w.r.t all possible realization of the uncertain variables $w_{i,t|k}^h$ and $w_{i,t|k}^c$ for all predicted time step $t\in\mathcal{T}$.

To extend the proposed formulation \eqref{opt1} to the energy management problem of smart thermal grids, we first need to introduce the notation, 
$\bs{x}_i := \{x_{i,t+1|k}\}_{t\in\mathcal{T}} \in \R^{6 N_h =: n_x}$, 
$\bs{u}_i:= \{u_{i,t|k}\}_{t\in\mathcal{T}}\in \R^{9 N_h =: n_u}$, 
$\bs{v}_i:=\{v_{i,t|k}\}_{t\in\mathcal{T}}\in \R^{2 N_h =: n_v}$, and 
$\bs{w}_i:=\{w_{i,t|k}\}_{t\in\mathcal{T}} \in \R^{2 N_h =: n_w}$. 
Given the initial value of the state $x_{i,k}\,$, one can eliminate the state variables from the dynamics \eqref{imb} of each agent $i$:
\begin{align}\label{}
\bs{x}_i = A_i x_{i,k} + B_i \bs{u}_i + C_{i} \bs{w}_i \,,
\end{align}
where the exact form of $A_i$, $B_i$ and $C_i$ matrices are omitted in the interest of space and can be found in \cite[Section 9.5]{borrelli2011predictive}.
We can now rewrite the total cost function over the prediction horizon in a more compact form as follows:
\begin{align*}
\cl{J}_i(\bs{x}_{i}, \bs{u}_i) = \bs{x}_i^\top \bs{Q}_i \bs{x}_i + \bs{u}_i^\top \bs{R}_i \bs{u}_i \ ,
\end{align*}
where $\bs{Q}_i$ and $\bs{R}_i$ are two block-diagonal matrices with $Q_i$ and $R_i$ on the diagonal for each agent $i$.
Note that the sum $\sum(\cdot)$ and the expectation $\Eb [\cdot]$ in the cost function \eqref{costopt1} are linear operators and thus, we can change their order without loss of generality.
Consider now the reformulation of \eqref{opt1} in a more compact form as follows:
\begin{subequations}
	\label{opt_agent}
	\begin{align}
	\underset{\bs{u}_{i}, \bs{v}_{i}}{\min} \quad & \cl{V}_i(\bs{x}_{i}, \bs{u}_i) = \Eb_{\bs{w}_i}\big[ \cl{J}_i(\bs{x}_{i}, \bs{u}_i) \big] \\
	\text{s.t.} \ \quad & E_i \bs{u}_i + F_i \bs{v}_i + P_i \leq 0 \,, \quad \forall\bs{w}_i \in \cl{W}_i\\
	\label{private_uncertain} &\Pb_{\bs{w}_i}\big[ A_i x_{i,k} + B_i \bs{u}_i + C_{i} \bs{w}_i \geq 0 \big] \geq 1-\varepsilon_i \,,
	\end{align}
\end{subequations}
where $E_i\,$, $F_i\,$, $P_i\,$ are matrices that are built by concatenating all constraints in \eqref{opt1}.
The index of $\Eb_{\bs{w}_i}, \Pb_{\bs{w}_i}$ denotes the dependency of the state trajectory $\bs{x}_{i}$ on the string of random scenarios $\bs{w}_i$ for each agent $i$.
The following technical assumption is adopted.  
\begin{assumption}
	The random variable $\bs{w}_i\,$ is defined on some probability space $(\cl{W}_i,\frk{B}(\cl{W}_i),\Pb_{\bs{w}_i})$, where $\cl{W}_i \subseteq \R^{n_w}$, $\frk{B}(\cdot)$ denotes a Borel $\sigma$-algebra, and $\Pb_{\bs{w}_i}$  is a probability measure defined over $\cl{W}_i$. 
\end{assumption}

It is worth to mention that for our study we only need a finite number of instances of $\bs{w}_i$, and we do not require the probability space $\cl{W}_i$ and the probability measure $\Pb_{\bs{w}_i}$ to be known explicitly.
The availability of a number of scenarios from the sample space $\cl{W}_i$ is enough which will become concrete in later parts of the paper.
Such samples can be for instance obtained from historical data.

The proposed optimization problem \eqref{opt_agent} is a finite-horizon, chance-constrained mixed-integer quadratic program, whose stages are coupled by the binaries \eqref{startup}, and dynamics of the imbalance error \eqref{probimb} for each agent $i$ at each sampling time $k$.
It is important to note that the proposed problem \eqref{opt_agent} is in general a non-convex problem and hard to solve.
In the following section, we will develop a tractable framework to obtain a probabilistically feasible solution for each agent $i$.
We refer to the proposed optimization problem \eqref{opt_agent} as a single agent problem, and whenever all agents solve this problem separately in a receding horizon fashion without any coupling constraints, it is referred to as the \textit{decoupled solution} (DS) in the subsequent parts.
We next extend the proposed single agent optimization problem \eqref{opt_agent} into a STGs setting.

\subsection{ATES in Smart Thermal Grids}

\begin{figure}
	\centering
	\includegraphics[width=\widthscale]{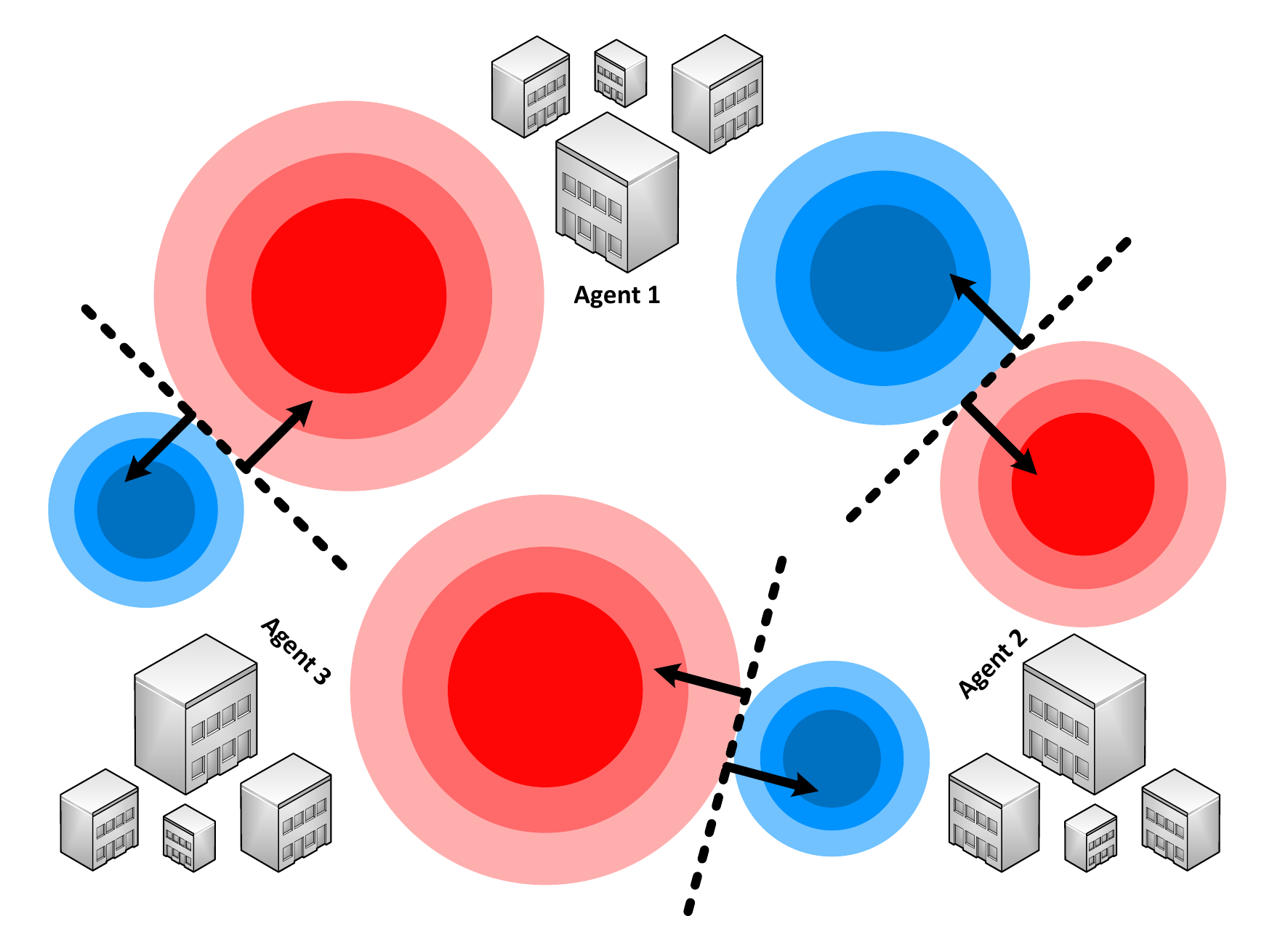}
	\caption{Three-agent ATES system in a STG. Each agent has a single ATES system which consists of a warm and a cold well. 
		Horizontal cross sections of warm and cold wells are shown with red and blue circles. 
		The black dashed lines represent the unwanted mutual interactions between ATES systems.}
	\label{fig_atessg}
\end{figure}

Consider a regional thermal grid consisting of $N$ agents with heterogeneous parameters as it was developed in the previous part. 
Such a STG setting however can lead to unwanted mutual interactions between ATES systems as it is illustrated in Fig.~\ref{fig_atessg}.
We therefore need to introduce a proper coupling constraint between neighboring agents that makes use of the following assumption. 

\begin{assumption}
	Each well of an ATES system is considered as a growing reservoir with respect to the horizontal axis (see black solid line in Fig.~\ref{fig_hcnets}).
	We therefore assume to have a cylindrical reservoir with a fixed height $\ell\, [\tx{m}]$ (filter screen length) and a growing radius $r_{\tx{a},k}^h \,, r_{\tx{a},k}^c\, [\tx{m}]$ (thermal radius) for each well of an ATES system.
\end{assumption}

Using the volume of stored water in each well of ATES system, one can determine the thermal radius using
\begin{align}\label{nonlinear_thermal_radios}
r_{\tx{a},k}^h = \left(\frac{c_{pw}\,\tx{V}_{\tx{a},k}^h}{c_{aq}\, \pi \ell}\right)^{0.5} \ , \ 	r_{\tx{a},k}^c = \left(\frac{c_{pw}\,\tx{V}_{\tx{a},k}^c}{c_{aq}\, \pi \ell}\right)^{0.5} \ , 
\end{align}
where $c_{aq} = (1-n_p) c_{\text{sand}} + n_p c_{pw}$ is the aquifer heat capacity.
$c_{\text{sand}}$ $[\tx{J} \text{kg}^{-1} \tx{K}^{-1}]$ relates to the sand specific heat capacity, and $n_p$ $[-]$ is the porosity of aquifer. 
Let us now denote the set of neighbors of agent $i$ by 
\begin{align*}
\cl{N}_i \subseteq \left\{ 1, 2, \cdots, N \right\} \backslash \{i\} \ .
\end{align*}
We impose a limitation on the thermal radius of warm well $r_{\tx{a},k}^h$ and cold well $r_{\tx{a},k}^c$ of ATES system in each agent $i$, based on the corresponding wells of its neighbor $j \in \cl{N}_i\,$:
\begin{align}\label{thermalradios}
(r_{\tx{a},k}^h)_i + (r_{\tx{a},k}^c)_j &\leq d_{ij} \ , \ j \in \cl{N}_i \ , 
\end{align}
where $d_{ij}$ is a given distance between agent $i$ and its neighbor $j\in \cl{N}_i$.
This constraint prevents overlapping between the growing domains of warm and cold wells of ATES systems in a STG setting.
Due to the nonlinear transformation in \eqref{nonlinear_thermal_radios}, we propose the following reformulation of this constraint to simplify the problem:
\begin{align}\label{volume}
(\tx{V}_{\tx{a},k}^h)_i + (\tx{V}_{\tx{a},k}^c)_j &\leq \tx{V}_{ij} - \bar\delta_{ij,k} \ ,
\end{align}
where $\tx{V}_{ij} = c_{aq}  \pi \ell \, (d_{ij})^2 / c_{pw}$ denotes the total volume of common resource pool between agent $i$ and its neighbor $j\in \cl{N}_i\,$. 
The variable $\bar\delta_{ij,k} = 2 c_{aq} \pi \ell \, (\bar{r}_{\tx{a},k}^h)_i \, (\bar{r}_{\tx{a},k}^c)_j / c_{pw}$ represents a time-varying parameter that captures the mismatch between the linear and nonlinear  constraint relations.
The following corollary is a direct result of the above reformulation. 

\begin{corollary}\label{cor_common_cons}
	If $(\bar{r}_{\tx{a},k}^h)_i$ and $(\bar{r}_{\tx{a},k}^c)_j$ represent the current thermal radius of warm and cold wells of ATES system in agent $i$ and $j$, respectively, then constraints \eqref{thermalradios} and \eqref{volume} are equivalent.
\end{corollary} 
The proof is provided in the Appendix.\QEDB

\begin{definition}\label{common_uncertainty_source}
	We define $\delta_{ij,k}$ to be a common uncertainty source between each agent $i$ and its neighboring agent $j\in \cl{N}_i$, using the following model:
	\begin{align}\label{common_unicertain_model}
	\delta_{ij,k} := \bar\delta_{ij,k} \,(1\pm 0.1 \, \zeta) \,,
	\end{align}
	where $\zeta$ is a random variable defined on some probability space, $\bar\delta_{ij,k}$ is constructed by using two given possible $(\bar{r}_{\tx{a},k}^h)_i \,, (\bar{r}_{\tx{a},k}^c)_j$ realizations that can be obtained using historical data in the DS framework.
	Since the mapping \eqref{common_unicertain_model} from $\zeta$ to $\delta_{ij,k}$ is measurable, one can view $\delta_{ij,k}$ as a random variable on the same probability space as $\zeta$.
\end{definition}

\subsection{Problem Formulation in Multi-Agent Network}

We now formulate the energy management problem for ATES systems in STGs as follows:
\begin{subequations}
	\label{opt_multiagent}
	\begin{align}
	\underset{\{\bs{u}_{i}, \bs{v}_{i}\}_{i=1}^{N}}{\min} \quad & \sum_{i=1}^{N} \ \cl{V}_i(\bs{x}_{i}, \bs{u}_i) \\
	\text{s.t.} \qquad\
	&\label{local_cons} E_i \bs{u}_i + F_i \bs{v}_i + P_i \leq 0 \,, \\
	&\label{w_uncer}\Pb_{\bs{w}_i}\left[ A_i x_{i,k} + B_i \bs{u}_i + C_{i} \bs{w}_i \geq 0 \right] \geq 1-\varepsilon_i \,, \\
	&\label{delta_uncer} \Pb_{\bs{\delta}_{ij}}  \left[ H_i \bs{x}_i + H_j \bs{x}_j \leq \bar{V}_{ij} - \bs{\delta}_{ij} \right] \geq 1 - \bar{\varepsilon}_{ij},  \\
	&\quad\forall \bs{w}_i \in \cl{W}_i  \,,\, \forall \bs{\delta}_{ij} \in \Delta_{ij}\,, \nonumber\\
	&\quad\forall j \in \cl{N}_i \,,\, \forall i \in \{1, 2, \cdots, N\} \,, &\nonumber
	\end{align}
\end{subequations}
where $H_i$, $H_j$ are coefficient matrices of appropriate dimensions, 
$\bar{V}_{ij}\in \R^{N_h}$ is the upper-bound on the total common resource pool, 
$\bs{\delta}_{ij}$ is a vector of common uncertainty variables,
and $\bar{\varepsilon}_{ij} \in (0,1)$ denotes the level of admissible coupling constraint violation for each agent $i$ and $\forall j\in \cl{N}_i$.
$\bar{V}_{ij}$ can be expressed as $\bar{V}_{ij} = \bs{1}^{N_h} \bigotimes V_{ij}$, using the Kronecker product.
Notice that the index of $\Pb_{\bs{\delta}_{ij}}$ denotes the dependency of the state trajectories on the string of random common scenarios 
$\bs{\delta}_{ij} = \{\delta_{ij,t|k}\}_{t\in\mathcal{T}}\subseteq \R^{N_h =: n_\delta}$.

\begin{assumption} The variable $\bs{\delta}_{ij}$ is considered to be a random vector on some probability space $(\Delta_{ij},\frk{B}(\Delta_{ij}),\Pb_{\bs{\delta}_{ij}})$, where $\Delta_{ij} \subseteq \R^{n_\delta}$, $\frk{B}(\cdot)$ denotes a Borel $\sigma$-algebra, and $\Pb_{\bs{\delta}_{ij}}$  is a probability measure defined over $\Delta_{ij}$. 
\end{assumption}

\begin{assumption}\label{independence}
	The variables $\bs{w}_{i}\in\R^{n_w}$ and $\bs{\delta}_{ij}\in\R^{n_\delta}$ are two vectors of independent random scenarios from two disjoint probability spaces $\cl{W}_{i}$ and $\Delta_{ij}$, respectively.
\end{assumption}

We refer to the proposed optimization problem \eqref{opt_multiagent} as a multi-agent network problem, and whenever the proposed problem \eqref{opt_multiagent} is solved in a receding horizon fashion, it is mentioned as the \textit{centralized solution} (CS) in the following parts.
The feasible set of \eqref{opt_multiagent} is in general non-convex and hard to determine explicitly due to the presence of chance constraints \eqref{w_uncer}, \eqref{delta_uncer}. 
In what follows, we will develop a tractable framework to obtain probabilistically feasible solutions for all agents.

\subsection{Move-Blocking Scheme}\label{multi}

The proposed system dynamics in \eqref{imb} for each agent $i$ consists of a BCC system dynamics \eqref{QB} with typically an hourly-based operation, and an ATES system \eqref{atesi} that is based on a seasonal variation of desired optimal operation.
This leads to a control problem that is sensitive w.r.t. the prediction horizon length, e.g., \eqref{opt_agent} and \eqref{opt_multiagent}.
Using a fixed prediction horizon length, e.g., least common multiple of these two systems, may turn out to be computationally prohibitive, however, also necessary in order to represent ATES interaction dynamics. 
We therefore aim to formulate a move-blocking strategy to reduce the number of control variables.

Consider $\cl{T} = \{k, k+1, \cdots, k+N_h-1\}$ to be the set of sampling time instances within the full prediction horizon, and 
$\cl{T}_u = \{\tau_1,\tau_2, \cdots,\tau_{T_u}\} \subseteq \cl{N}_h$ to be the set of sampling instances at which the control input is updated with $T_u = |\cl{T}_u |$.
We introduce a new vector of multi-rate decision variables $\tilde{\bs{u}}_i \in \R^{N_u T_u}$ which are related to the original ones by: 
\begin{align}\label{multi_trans}
\bs{u}_i = \bs{\Psi}\tilde{\bs{u}}_i \ ,
\end{align}
where  
$\bs{\Psi} = \left[\begin{matrix} \Psi_{1} & \Psi_{2} & \cdots & \Psi_{T_u} \end{matrix}\right] \in \R^{N_u N_h \times N_u T_u}$ is a linear mapping matrix.
For all $m\in\{1, \cdots, T_u\}$, we construct 
\begin{align}
\Psi_{m} = \left[\begin{matrix} \psi_{1,m}^\top & \psi_{2,m}^\top & \cdots & \psi_{N_h,m}^\top \end{matrix}\right]^\top \in \R^{N_u N_h \times N_u} \ ,
\end{align}
where $\psi_{l,m} \in \R^{N_u \times N_u}$ for all $l \in \{1, 2, \cdots, N_h\}$ is defined as
\begin{align}\label{multi_rate}
\psi_{l,m} = 
\begin{cases}
\bs{1} & \text{if } k+l-1 = \tau_m	\\
0 				 & \text{otherwise}
\end{cases}\,,
\end{align} 
where $\bs{1} \in \R^{N_u \times N_u}$ represents an identity matrix.

We reformulate the optimization problem \eqref{opt_multiagent} using the proposed move-blocking scheme \eqref{multi_rate}, and whenever the reformulation of \eqref{opt_multiagent} is solved in a receding horizon fashion, it is referred to as the \textit{move-blocking centralized solution} (MCS).

\section{Computationally Tractable Framework}\label{sec_framework}

In this section, we provide a framework to approximately solve the chance-constrained optimization problem \eqref{opt_multiagent}, which is in general difficult to solve.
To this end, we employ a data-driven approach to approximate the chance constraints using some available samples of uncertainties.
We first extract at random some instances of the uncertainties (scenarios), where the scenarios are independent and identically distributed (i.i.d.) and, then, find the optimal solution of the problem with only the constraints associated with the extracted scenarios. 

An important requirement of such a technique is to have a convex problem w.r.t the decision variables, which is not the case in our formulation \eqref{opt_multiagent}.
To tackle such a mixed-integer chance-constrained problem, one can use a worst-case mixed-integer reformulation technique as it was initially introduced in \cite{esfahani2015performance}.
Due to the large-scale network problem \eqref{opt_multiagent}, such a reformulation leads to enormous cost of computation and it is indeed an intractable approach.
Following the so-called robust randomized technique \cite{margellos2014road}, the reformulation is done in a way to provide a feasible solution for all scenarios of the uncertainty realizations in a probabilistic sense. 

The idea of robust randomized approach is the following. 
An auxiliary chance-constrained optimization problem is first formulated to determine a probabilistic bounded set of random variables.
This yields a bounded set of uncertainty that is a subset of the uncertainty space and contains a portion of the probability mass of the uncertainty with high confidence level. 
Then, a robust version of the initial problem subject to the uncertainty confined in the obtained set is solved.
We here extend this framework in order to be able to handle a problem with multiple chance constraints based on the idea of the robust randomized approach\cite{margellos2014road}.

Consider $\bs{y}_i = (\bs{u}_i, \bs{v}_i) \in \R^{(n_u + n_v) = n_y}$,  $\bs{y} = \col(\bs{y}_i)_{i=1}^{N}$, where $\col(\cdot)$ is an operator to stack elements.
Define 
$\bs{w} = \col(\bs{w}_i)_{i=1}^{N} \subseteq \cl{W}$ to be the private uncertainty sources for a network of agents, 
$\bs{\delta}_i = \col(\bs{\delta}_j)_{j\in \cl{N}_i}^{} \subseteq \Delta_i$ to be the common uncertainty sources for each agent, and 
$\bs{\delta} = \col(\bs{\delta}_i)_{i=1}^{N} \subseteq \Delta$ to be the common uncertainty sources for a multi-agent network, where
\begin{align*}
\cl{W} := \prod_{i=1}^{N} \, \cl{W}_i \ , \
\Delta_i := \prod_{j\in \cl{N}_i}^{} \, \Delta_{ij} \ , \
\Delta := \prod_{i=1}^{N} \, \Delta_i \ . 
\end{align*}
Consider now the proposed optimization problem in \eqref{opt_multiagent} in a more compact format:
\begin{subequations}
	\label{comp_opt_multiagent}
	\begin{align}
	\underset{\bs{y}}{\text{min}} \quad & \sum_{i=1}^{N}  \cl{V}_i(\bs{x}_{i}, \bs{u}_i) \\
	\text{s.t.} \quad 
	&\label{private_uncer} \Pb_{\bs{w}}  \left[ \bs{y} \in \prod_{i=1}^{N} \cl{Y}_{i}(\bs{w}_{i}) \right] \geq 1-\varepsilon \, , \ \forall \bs{w} \in \cl{W} \\
	&\label{common_uncer} \Pb_{\bs{\delta}}  \left[ \bs{y} \in \prod_{i=1}^{N} \bigcap_{j\in \cl{N}_i}^{} \breve{\cl{Y}}_{ij}(\bs{\delta}_{ij}) \right] \geq 1 - \bar{\varepsilon} \, ,  \ \forall \bs{\delta} \in \Delta
	\end{align}
\end{subequations}
where $\varepsilon := \sum_{i=1}^{N} \varepsilon_i \in (0,1)$, $\bar{\varepsilon} := \sum_{i=1}^{N} \sum_{j\in \cl{N}_i}^{} \bar{\varepsilon}_{ij} \in (0,1)$.
$\cl{Y}_{i}(\bs{w}_{i})\in\R^{n_y}$ and $\cl{Y}_{ij}(\bs{\delta}_{ij})\in\R^{n_y}$ are defined\footnote{Both sets have a dependency on the initial value of the state $x_{i,k}$ for each agent $i$ at each sampling time $k$. Given $x_{i,k}$, we here highlight the dependency of these sets on the uncertainties $\bs{w}_{i}$ and $\bs{\delta}_{ij}$ for each agent $i$ at each sampling time $k$.} by
\begin{align*}
\cl{Y}_{i}(\bs{w}_{i}) &:= \big\{ \bs{y}_i\in\R^{n_y} \, : \, E_i \bs{u}_i + F_i \bs{v}_i + P_i \leq 0  \,, \\  
&\qquad\qquad\qquad\quad A_i x_{i,k} + B_i \bs{u}_i + C_{i} \bs{w}_i \geq 0 \big\} \,, \\
\cl{Y}_{ij}(\bs{\delta}_{ij}) &:= \big\{ \left(\bs{y}_i, \bs{y}_j\right) \in\R^{2 n_y} \, : \, H_i \bs{x}_i + H_j \bs{x}_j \leq \bar{V}_{ij} - \bs{\delta}_{ij}  \big\}.
\end{align*}
It is important to note that $\breve{\cl{Y}}_{ij}(\bs{\delta}_{ij})\in\R^{2 n_y N_i}$ represents the cylindrical extension\footnote{Cylindrical extension replicates the membership degrees from the existing dimensions into the new dimensions \cite[\S 4]{zadeh1997toward}.} of $\cl{Y}_{ij}(\bs{\delta}_{ij})$.
In the subsequent parts, we refer to the constraint \eqref{private_uncer} as the agents' private chance constraints, and to the constraint \eqref{common_uncer} as the agents' common chance constraints.
The proposed formulation \eqref{comp_opt_multiagent} is a mixed-integer quadratic optimization problem with multiple chance constraints, due to the binary variables $\{\bs{v}_i\}_{i=1}^{N}$ and the chance constraints \eqref{private_uncer}, \eqref{common_uncer}.
The index of $\Pb_{\bs{w}}$ and $\Pb_{\bs{\delta}}$ denote the dependency on the string of random scenarios $\bs{w} \in \cl{W}$ and $\bs{\delta} \in \Delta$, respectively.

Building upon our previous work in \cite{rostampour2016robust}, we extend the so-called robust randomized approach in \cite{margellos2014road,margellos2013stochastic}
to be able to handle a problem with multiple chance constraints.  
Problem \eqref{comp_opt_multiagent} is a stochastic program with multiple chance constraints, where $\Pb_{\bs{w}}$ and $\Pb_{\bs{\delta}}$ denote two different probability measures for private and common uncertainty sources, respectively. 
In summary, one can reformulate the chance constraints in \eqref{comp_opt_multiagent} using a worst-case chance constraint defined by
\begin{align}\label{mcp_scp}
\max_{\eta \in \cl{N}_{\text{MCP}}} \Pb \left[ f_\eta(\bs{y}, \cdot) \right] \geq 1 - \tilde{\varepsilon} \ ,
\end{align}
where $\tilde{\varepsilon} = \min_{\eta \in \cl{N}_{\text{MCP}}} \{\varepsilon_\eta\}$, 
$f_\eta(\bs{y}, \cdot)$ denotes the $\eta$-th chance constraint function, and $\cl{N}_{\text{MCP}}$ is the set of indices of chance constraint functions formulated in \eqref{comp_opt_multiagent}.
However, this procedure clearly leads to a considerable amount of conservatism, due to the fact that it requires the solution to satisfy all constraints with the highest probability $1 - \tilde{\varepsilon}$.
We instead employ the robust randomized approach for each chance constraint function 
$f_k(\bs{y}, \cdot), k \in \cl{N}_{\text{MCP}}$, separately.
Our framework is also related to, albeit different from the approach of \cite{schildbach2013randomized}, since the feasible set in \eqref{comp_opt_multiagent} is non-convex.  
Moreover, the problem formulation  in \cite{schildbach2013randomized} consists of an objective function with multiple chance constraints, in which the terms in objective and constraints are univariate w.r.t. the decision variables.
In contrast, the objective function in our problem formulation \eqref{comp_opt_multiagent} consists of separable additive components and constraint functions are also separable w.r.t. \eqref{private_uncer}, \eqref{common_uncer} between each agent $i =1, \cdots, N$ and $\forall j\in \cl{N}_i$.

Define $\cl{B}_{i}$, $\bar{\cl{B}}_{ij}$  to be two bounded sets of private uncertainty source and a bounded set of common uncertainty source for each agent $i$, respectively.
$\cl{B}_{i}$, $\bar{\cl{B}}_{ij}$ are assumed to be axis-aligned hyper-rectangular sets \cite[Proposition 1]{margellos2014road}.
This is not restrictive and any convex set with convex volume could have been chosen instead as in \cite{rostampour2017set}.
We parametrize $\cl{B}_{i} (\bs{\gamma}) :=[\bs{\ogam},\bs{\ugam}]$ by $\bs{\gamma}=(\bs{\ogam},\bs{\ugam})\in\R^{2n_w}$, and $\bar{\cl{B}}_{ij} (\bs{\lambda}) :=[\bs{\olam},\bs{\ulam}]$ by $\bs{\lambda}=(\bs{\olam},\bs{\ulam})\in\R^{2n_\delta}$, and 
consider the following chance-constrained optimization problem:
\begin{subequations}\label{PB}
	\begin{align}
	&\begin{cases}\label{PBi}
	\min\limits_{\bs{\gamma}} & \left\|\bs{\ogam} - \bs{\ugam}\right\|_1  \\ 
	\ \text{s.t.}  &\Pb \left\{\, \bs{w}_{i} \in \cl{W}_i \ \big| \ \bs{w}_{i} \in [\bs{\ugam}, \bs{\ogam}] \, \right\} \geq 1- \varepsilon_i
	\end{cases} \,, \\
	&\begin{cases}\label{PBij}
	\min\limits_{\bs{\lambda}_{}} & \left\|\bs{\olam}_{} - \bs{\ulam}_{}\right\|_1  \\ 
	\ \text{s.t.}  &\Pb \left\{\, \bs{\delta}_{ij} \in \Delta_{ij} \ \big| \ \bs{\delta}_{ij} \in [\bs{\ulam}_{}, \bs{\olam}_{}] \, \right\} \geq 1- \bar\varepsilon_{ij}
	\end{cases}\,.
	\end{align}
\end{subequations}
Following the so-called scenario approach in \cite{calafiore2006scenario}, one can determine the number of required uncertainty scenarios to formulate a tractable problem, using $N_s = \frac{2}{\epsilon} (\xi + \ln \frac{1}{\nu})$, where $\xi$ is the dimension of decision vector, $\epsilon$, $\nu$ are the level of violation, and the confidence level, respectively.
We determine $N_{s_i}$ by substituting $\xi = 2n_w$, $\epsilon = \varepsilon_i$, $\nu = \beta_i$, and determine $\bar{N}_{s_{ij}}$ by substituting 
$\xi = 2n_\delta$, $\epsilon = \bar\varepsilon_{ij}$, $\nu = \bar\beta_{ij}$, for all agent $i \in\{1, 2, \cdots, N\}$.

We next define $\cl{S}_i  = \{\bs{w}^{(1)}_{i}, \cdots, \bs{w}^{(N_{s_i})}_{i}\} \subset \cl{W}_{i}$, 
$\bar{\cl{S}}_{ij}  = \{\bs{\delta}^{(1)}_{ij}, \cdots, \bs{\delta}^{(\bar{N}_{s_{ij}})}_{ij}\} \subset \Delta_{ij}$ 
and formulate a tractable version of  \eqref{PBi} and \eqref{PBij} by 
\begin{subequations}\label{SPB}
	\begin{align}
	&\begin{cases}\label{SPBi}
	\min\limits_{\bs{\gamma}} & \left\|\bs{\ogam} - \bs{\ugam}\right\|_1  \\ 
	\ \text{s.t.}  & \bs{w}_{i} \in [\bs{\ugam}, \bs{\ogam}] \quad , \quad \bs{w}_{i} \in \cl{S}_i 
	\end{cases} \ ,  \\
	&\begin{cases}\label{SPBij}
	\min\limits_{\bs{\lambda}_{}} & \left\|\bs{\olam}_{} - \bs{\ulam}_{}\right\|_1  \\ 
	\ \text{s.t.}  & \bs{\delta}_{ij} \in [\bs{\ulam}_{}, \bs{\olam}_{}] \quad , \quad  \bs{\delta}_{ij} \in \bar{\cl{S}}_{ij} 
	\end{cases} \ .
	\end{align}
\end{subequations}
The optimal solutions $(\bs{\gamma}^*\,, \bs{\lambda}^*)$  of the proposed tractable problem are probabilistically feasible for the chance-constrained problems, \cite[Theorem 1]{campi2008exact}.
Moreover, $\bs{\gamma}^*$, and $\bs{\lambda}^*$ also characterize our desired probabilistic  bounded sets $\cl{B}_{i}^*$ and $\bar{\cl{B}}_{ij}^*$, respectively. 
Note that $\cl{S}_i$ and $\bar{\cl{S}}_{ij}$ are two collections of random scenarios that are i.i.d.

After determining $\cl{B}_{i}^*$ and $\bar{\cl{B}}_{ij}^*$ for all agents $i \in\{1, \cdots, N\}$, we are now able to reformulate the robust counterpart of the original problem \eqref{comp_opt_multiagent} via:
\begin{subequations}
	\label{robust_comp_opt_multiagent}
	\begin{align}
	\underset{\bs{y}}{\min} \quad & \sum\limits\nolimits_{i=1}^{N}  \cl{V}_i(\bs{x}_{i}, \bs{u}_i) \\
	\ \text{s.t.} \quad 
	&\bs{y} \in \, \prod\limits_{i=1}^{N} \, \bigcap\limits_{\bs{w}_i\in \{\cl{B}_{i}^* \bigcap \cl{W}_i\}} \, \cl{Y}_{i}(\bs{w}_{i}) \ , \\
	&\bs{y} \in \, \prod\limits_{i=1}^{N} \, \bigcap\limits_{j\in \cl{N}_i}^{} \, \bigcap\limits_{\bs{\delta}_{ij}\in \{\bar{\cl{B}}_{ij}^* \bigcap \Delta_{ij}\}} \, \breve{\cl{Y}}_{ij}(\bs{\delta}_{ij})  \ . \label{couple_pool}
	\end{align}
\end{subequations}
Note that the aforementioned problem is not a randomized program, and instead, the constraints have to be satisfied for all values of the private uncertainty in $\{\cl{B}_{i}^* \bigcap \cl{W}_i\}$, and common uncertainty in $\{\bar{\cl{B}}_{ij}^* \bigcap \Delta_{ij}\}$.
The proposed problem \eqref{robust_comp_opt_multiagent} is a robust mixed-integer quadratic program.
In \cite{bertsimas2006tractable}, it was shown that robust problems are tractable \cite[Proposition 1]{rostampour2016robust}, and remain in the same class as the original problems, e.g., robust mixed-integer programs remain mixed-integer programs, for a certain class of uncertainty sets, such as in our problem \eqref{robust_comp_opt_multiagent}, the uncertainty is bounded in a convex set. 
The following theorem quantifies the robustness of solution obtained by \eqref{robust_comp_opt_multiagent} w.r.t. the initial problem \eqref{comp_opt_multiagent}.

\begin{theorem}\label{thm_common}
	Let 
	$\varepsilon_{i}\,,\,\bar\varepsilon_{ij} \in (0,1)$ and 
	$\beta_{i}\,,\,\bar\beta_{ij} \in (0,1)$ 
	for all $ j\in\cl{N}_i$\,, for each $ i \in\{1, \cdots, N\}$ be chosen such that
	$\varepsilon = \sum_{i=1}^{N} \varepsilon_i\in (0,1)$\,,\,$\beta = \sum_{i=1}^{N} \beta_i\in (0,1)$\,,
	$\bar\varepsilon_i = \sum_{j\in\cl{N}_i}^{} \bar\varepsilon_{ij}\in (0,1)$\,,\,$\bar\beta_i = \sum_{j\in\cl{N}_i}^{} \bar\beta_{ij}\in (0,1)$ and
	$\bar\varepsilon = \sum_{i=1}^{N} \bar\varepsilon_i\in (0,1)$\,,\,$\bar\beta = \sum_{i=1}^{N} \bar\beta_i\in (0,1)$.
	Determine $\cl{B}_{i}^*$ and $\bar{\cl{B}}_{ij}^*$ by constructing ${\cl{S}}_{i}$\,,\,$\bar{\cl{S}}_{ij}$ and solving \eqref{SPB} for all $j\in\cl{N}_i$\,, for each $i \in\{1, \cdots, N\}$.
	If $\bs{y}^*_{s}$ is a feasible solution of the problem \eqref{robust_comp_opt_multiagent}, then $\bs{y}^*_{s}$ is also a feasible solution for the chance constraints \eqref{private_uncer} and \eqref{common_uncer}, with the confidence levels of $1-\beta$ and $1-\bar\beta$, respectively.
\end{theorem}
The proof is provided in the Appendix.\QEDB

The interpretation of \Cref{thm_common} is as follows. 
The obtained solutions via \eqref{robust_comp_opt_multiagent}
for all agents $i=1,\cdots,N$ have feasibility guarantees with $1-\varepsilon_i$ and $1-\bar{\varepsilon}_{ij}$ probabilities for the private and common uncertain sources $\bs{w}_i$ and $\bs{\delta}_{ij}$ with high confidence levels $\beta_i$ and $\bar{\beta}_{ij}$, respectively. 
To keep the robustness level of the solutions for the whole network problem, these choices have to follow a certain design rule.
It is important to mention that in order to maintain the violation level for the whole network the violation level of individual agent needs to decrease which may lead to very conservative results for each agent, since the number of required samples needs to increase in the proposed formulation \eqref{SPB}.

\begin{remark}
	We also approximate the objective function empirically for each agent $i$ following the approach in \cite{rostampour2015stochastic}.
	$\Eb_{\bs{w}_i}[\cl{J}_i(\cdot)]$ can be approximated by averaging the value of its argument for some number of different scenarios, which plays a tuning parameter role.  
	Using $N_{s_{i}^0}$ as the tuning parameter, consider $N_{s_i^0}$ number of different scenarios of $\bs{w}_{i}$ to build $\cl{S}_i^0  = \{\bs{w}^{(1)}_{i}, \cdots, \bs{w}^{(N_{s_i^0})}_{i}\} \subset \cl{W}_{i}$ for each agent $i = 1, \cdots, N$. 
	Then approximate the cost function empirically as follows: 
	\begin{align*}
	\sum\nolimits_{i=1}^{N} \cl{V}_i(\bs{x}_{i}(\bs{w}_i),\bs{u}_{i}) =  \sum\nolimits_{i=1}^{N} \Eb_{\bs{w}_i\in \cl{W}_i}\left[\cl{J}_i(\bs{x}_{i}(\bs{w}_i), \bs{u}_{i})\right] 
	\approx  \sum\nolimits_{i=1}^{N} \frac{1}{N_{s_i^0}} \sum\nolimits_{\bs{w}_i\in \cl{S}_i^0} \cl{J}_i(\bs{x}_{i}(\bs{w}_i), \bs{u}_{i}) \,.
	\end{align*}
	It is worth mentioning that one can employ scenario removal algorithms to improve the objective value, leading to a tradeoff between feasibility and optimality, see e.g., \cite{esfahani2015performance, campi2011sampling}.
	
\end{remark}

\begin{remark}
	A tractable decoupled solution (DS) formulation for \eqref{opt_agent} can be achieved by removing the robust coupling constraint \eqref{couple_pool} from \eqref{robust_comp_opt_multiagent}.
	Since there is no longer a coupling constraint, each agent $i$ can therefore solve its problem independently. 
\end{remark}

The solution of \eqref{robust_comp_opt_multiagent} is the optimal input sequence $\{u_{i,k|k}^{*},  v_{i,k|k}^{*}, \cdots, u_{i,k+N_{h}-1|k}^{*}, v_{i,k+N_{h}-1|k}^{*}\}_{i=1}^{N}$.
Based on an MPC paradigm, the current input at time step $k$ is implemented in the system dynamics \eqref{imb} using the first element of optimal solutions as $\{u_{i,k},  v_{i,k}\}_{i=1}^{N} := \{u_{i,k|k}^{*},  v_{i,k|k}^{*}\}_{i=1}^{N}$ and we proceed in a receding horizon fashion. 
This means \eqref{robust_comp_opt_multiagent} is solved at each time step $k$ by using the current measurement of the state $\{x_{i,k}\}_{i=1}^{N}$. 
It is important to highlight that the feasibility guarantees in \Cref{thm_common} are independent from the sampling rate of the real continuous-time system.
It is however very important to have a discrete-time system model that can predict the real system behavior as precisely as possible.	
Once such a suitable discrete-time system model is developed, one can use our proposed tractable frameworks (DS, CS, and MCS), and instead of analyzing the closed-loop asymptotic behavior, achieve the fulfillment of multiple chance constraints from an optimization point of view and have a-priori probabilistic feasibility guarantees via \Cref{thm_common}.

\section{Numerical Study}\label{sec_sim}

In this section, we present a simulated case study for a three-agent ATES system in a STG, as it is shown in Fig.~\ref{fig_atessg}.
We determine the thermal energy demands of three buildings, that had been equipped with ATES systems, modeled using realistic parameters and the actual registered weather data in the city center of Utrecht, The Netherlands, where these buildings are located. 
We refer interested readers to \cite[Appendix A]{wicak2016distributed} for the complete detailed parameters of this case study.

\subsection{Simulation Setup}

We simulate three problem formulations, namely: DS (decoupled solution), CS (centralized solution), and MCS (move-blocking centralized solution), using the proposed tractable framework \eqref{robust_comp_opt_multiagent}.
The simulation time is one year from June 2010 to June 2011 with hourly-based sampling time.
The prediction horizon for DS and CS is a day-ahead (24 hours), whereas for MCS is a whole season (3 months). 
The multi-rate control actions in MCS are considered to be hourly-based during first day, daily-based in the first week, weekly-based within the first month, and monthly-based for the rest of the season.
We also simulate a deterministic DS (DDS) for comparison purposes, where the uncertain elements $(\bs{w}_i)$ are fixed to their forecast value for each agent $i = 1, 2, 3$.
In order to generate scenarios from the private uncertainty sources, we use a discrete normal stochastic process, where the thermal energy demand of each building varies within 10\% of its actual value at each sampling time. A similar technique is used for the common uncertainty sources.
The simulation environment was MATLAB with YALMIP as the interface \cite{lofberg2004yalmip} and Gurobi as a solver.

\subsection{Simulation Results}

\begin{figure}
	\centering
	\includegraphics[width=\widthscale]{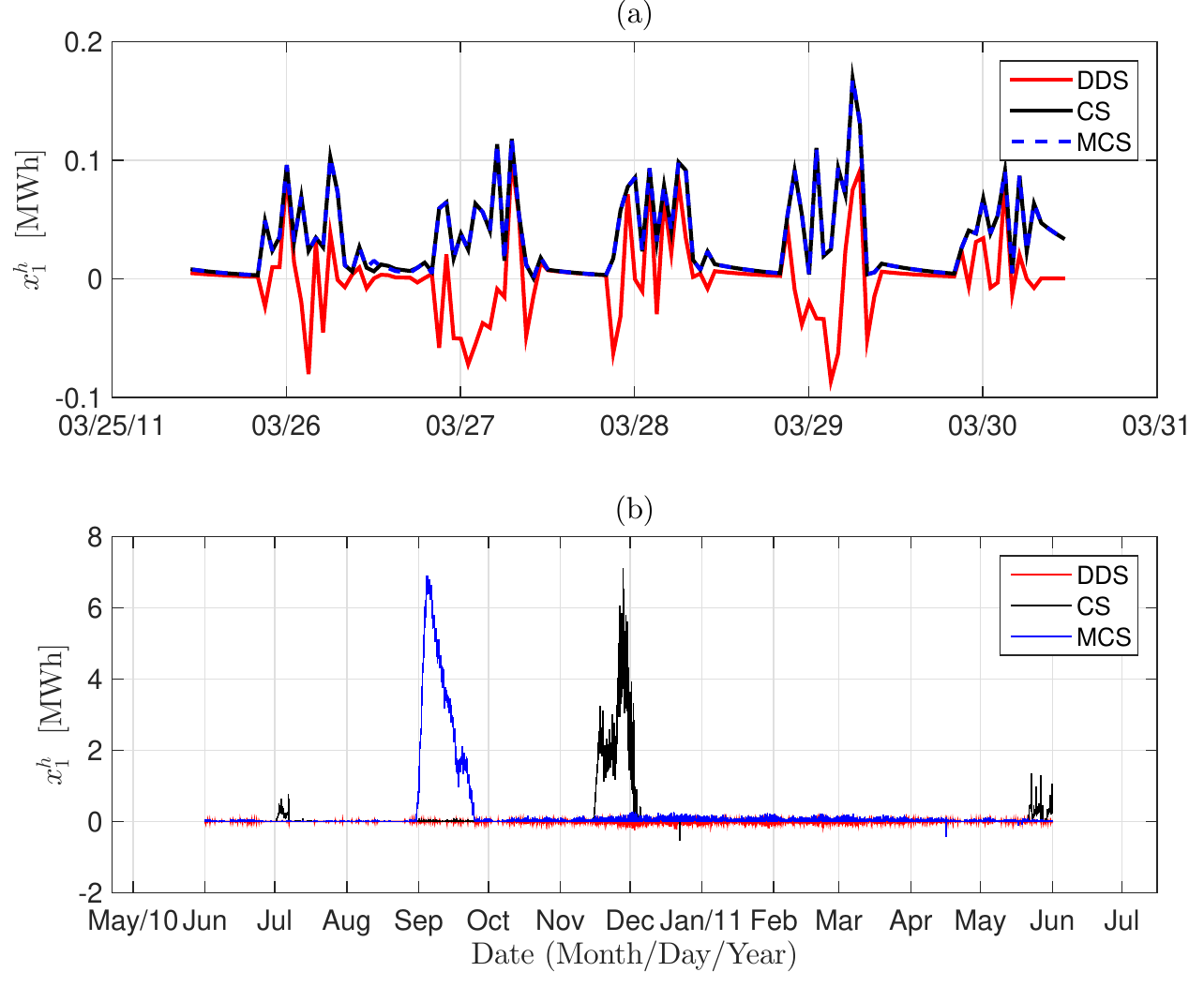}
	\caption{A-posteriori feasibility validation of the obtained results via DDS, CS,
		and MCS formulations for the imbalance error dynamics in the first building of the three-agent ATES-STG example. Fig.3(a)  focuses on a randomly chosen five-day period to allow a better comparison, whereas Fig.3(b) presents the complete one year results.}
	\label{fig_res1}
\end{figure}

\begin{figure}
	\centering
	\includegraphics[width=\widthscale]{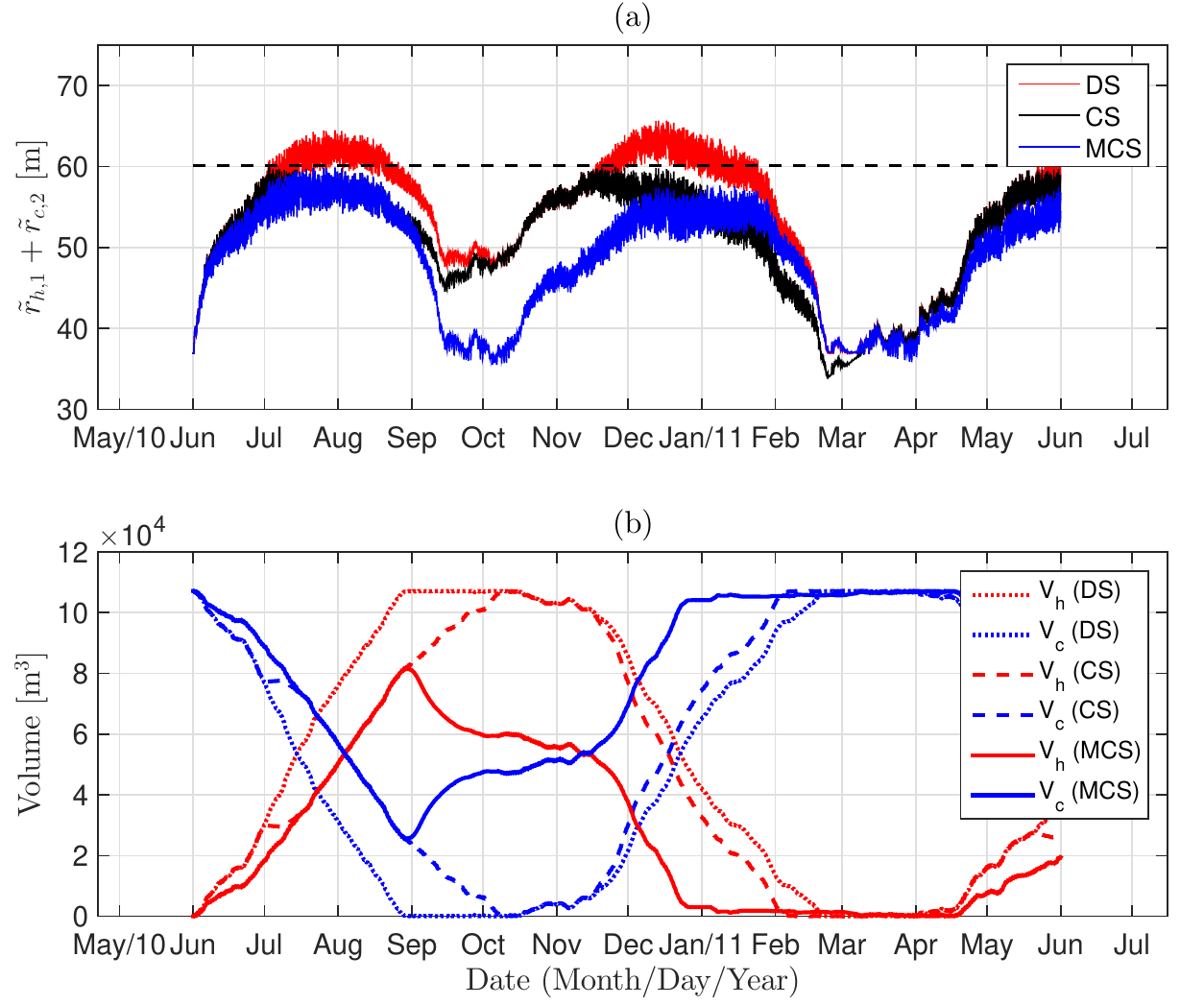}
	\caption{A-posteriori feasibility validation of the
		obtained results via DS, CS, and MCS formulations for the common coupling
		constraint between the first and second building of the three-agent ATES-STG
		example is presented in Fig.4(a). Fig.4(b) shows the ATES system state trajectories
		(volume of the stored water in the warm and cold wells of the first building)
		in the three-agent ATES-STG example.}
	\label{fig_res2}
\end{figure}

Fig.~\ref{fig_res1} and Fig.~\ref{fig_res2} (a) depict a-posteriori feasibility validation of the private chance constraint of agent 1 and the common chance constraint between agent 1 and agent 2.
It is important to note that the results obtained for the other two buildings are very similar, and therefore we focus on the results of the first building (agent 1).
To illustrate the functionality of our proposed framework to deal with the private chance constraint, in Fig.~\ref{fig_res1}, we present the a-posteriori feasibility validation of the obtained results via DDS, CS, and MCS formulations.
Fig.~\ref{fig_res1} (a) shows the obtained results for the last five days in March 2011, and Fig.~\ref{fig_res1} (b) shows the results for one year simulation from June 2010 until June 2011.
In Fig.~\ref{fig_res1} the "red" color denotes the solution of DDS, "black" color shows the solution of CS, and "blue" presents the solution of MCS.

Fig.~\ref{fig_res1} (a) focuses on a randomly chosen five-day period to allow a better comparison between the results of DDS, CS, and MCS. 
It is clearly shown that the obtained results via CS and MCS, provide a feasible (nonnegative) trajectory of the thermal energy imbalance error during heating mode, whereas the solution of DDS, leads to some violations throughout the simulation time.
Notice that all three proposed approaches, namely DS, CS, and MCS, achieved the feasibility of the private chance constraint in a probabilistic sense as it is guaranteed in Theorem \ref{thm_common}.
We present the results obtained via DDS to highlight such an achievement, whereas the results obtained via DS is omitted to demonstrate the other achievements.

In Fig.~\ref{fig_res1} (b), the complete one year results of DDS, CS, and MCS are shown. 
Two important observations are as follows: the obtained results of CS and MCS have very small number of violations, much less than our desired level  of violations, throughout the simulation time. 
This yields a less conservative approach compared to the classical robust control approach (see \cite[Ch.14]{borrelli2011predictive}).
As the second observation, in the results of CS and MCS one can see some instances of a large non-zero imbalance error, which is expected:
By taking into account the coupling constraints between agents, the solutions of agents are going to extract the stored thermal energy from their ATES systems to prevent the mutual interactions between their ATES systems as in Fig.~\ref{fig_res2} (a). 
Interestingly, the results of MCS show that agent 1 starts to extract the stored thermal energy from its ATES system sooner due to its longer prediction horizon, compared to CS.

Fig.~\ref{fig_res2} (a) shows the evaluation of our proposed reformulation for the coupling constraint in \eqref{volume} together with the a-posteriori feasibility validation of the common chance constraint between agent 1 and agent 2. 
We plot the obtained $\tilde{r}_{h,1} + \tilde{r}_{c,2}$ using DS, CS, and MCS formulations. 
As it is clearly shown DS results are violating the coupling constraint which leads to overlap between the stored water in warm well of ATES system in agent 1 and the stored water in cold well of ATES system in agent 2.
This is due to the fact that there are no coupling constraints in the DS framework and each agent works without any information from neighboring agents.
It is important to highlight that the results obtained via DDS and DS are the same in terms of the ATES system dynamical behavior.
This is due to the fact that the cost parameter associated with the ATES system pump is the same in both DDS and DS formulations, and thus ATES systems participate in the agent energy management in the same way, regardless of the private chance constraints.   
We also present the evolution of the stored water volume in each well of the ATES system for agent 1 using the obtained results via DS, CS, and MCS formulations in Fig.~\ref{fig_res2} (b) to illustrate the impact of the different formulations.

It is worth to mention that Fig.~\ref{fig_res1} and Fig.~\ref{fig_res2} illustrate all main contributions: 
1) having a probabilistically feasible solution for each agent w.r.t. the private uncertainty sources as it is encoded via \eqref{private_uncer},
2) respecting the common resource pool between neighboring agents in STGs as it is formulated in \eqref{common_uncer} 
(the first and second outcomes are the direct results of our theoretical guarantee in Theorem \ref{thm_common}), and
3) prediction using a longer horizon yields an anticipatory control decision that improves the operation of an ATES system. 
This is a direct consequence of our proposed move-blocking scheme in \eqref{multi_rate}.

\begin{figure}
	\centering
	\includegraphics[width=\widthscale]{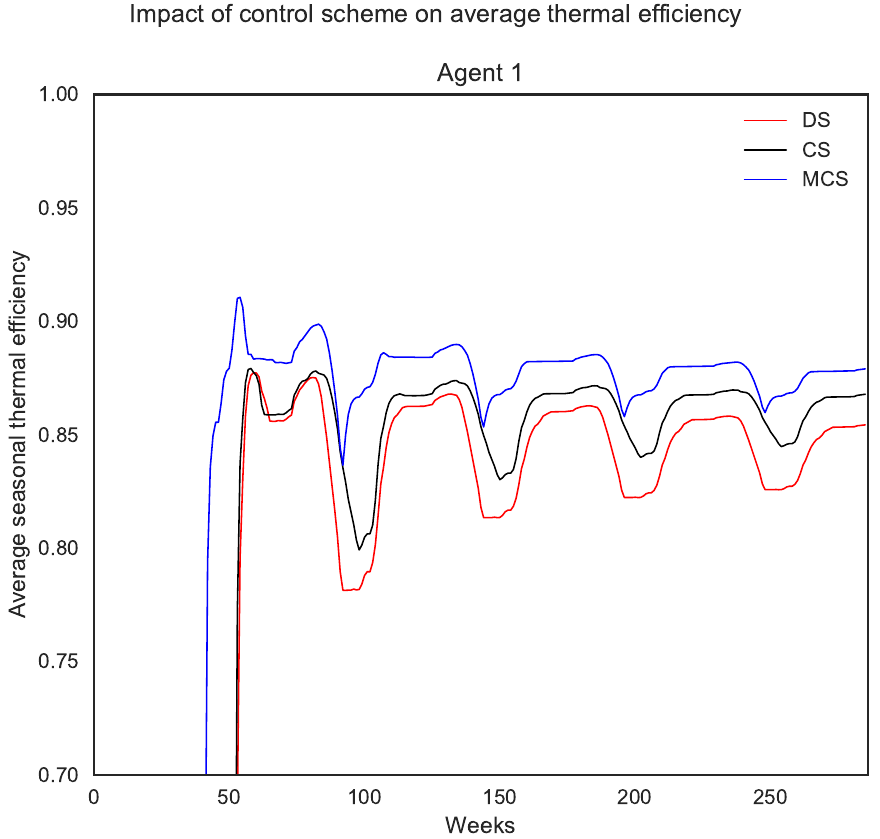}
	\caption{Impact of DS and CS on average thermal efficiency.}
	\label{fig_marc}
\end{figure}

Fig.~\ref{fig_marc} summarizes the results in terms of average thermal efficiency that we obtained by integrating our control strategy, DS and CS, into Python to build a live-link with MODFLOW, 
a more realistic aquifer
simulation environment\footnote{MODFLOW is a modular hydrologic model, and it is considered an international standard for aquifer simulation and predicting groundwater conditions and interactions.} \cite{harbaugh2005modflow}.
Fig.~\ref{fig_marc} is presented to highlight the impact of considering the proposed coupling constraints, as it is formulated in \eqref{common_uncer}, versus the decoupled setting.
The impact of our control strategy, DS (red) and CS (blue), on average thermal energy efficiency \cite{bloemendal2018methods} in each building illustrates that we can store and retrieve the same amount of thermal energy in ATES systems, in a more efficient way due to information exchange between the agents to prevent the mutual interactions between wells using the results of MCS and CS compared to DS.

\section{Conclusions and Future Work}\label{final}

This paper proposed a stochastic MPC framework for an energy management problem in STGs consisting of ATES systems integrated into BCC systems.  
We developed a large-scale stochastic hybrid model to capture thermal energy imbalance errors in an ATES-STG.
In such a framework, we formalized two important practical concerns, namely: 1) the balance between extraction and injection of energy from and into the aquifers within a certain period of time; 
2) the unwanted mutual interaction between ATES systems in STGs.
Using our developed model, we formulated a finite-horizon mixed-integer quadratic optimization problem with multiple chance constraints. To solve such a problem, we proposed a tractable formulation based on the so-called robust randomized approach.
In particular, we extended this approach to handle a problem with multiple chance constraints.
We simulated our proposed framework using a three-agent ATES-STG example which confirmed the expected performance improvements.  
Our current work focuses on two main directions:

a) From an application point of view, refining the proposed model of ATES system \eqref{ates_dyn} to be able to predict situations where the wells are completely depleted, or the new ones are installed. Such a situation may happen in reality where water extraction is continued with the aquifer ambient temperature. In \cite{rostampour2017aquifer}, we developed such a model and currently the possibility of integration of such a development into STGs is under investigation.

b) From a theoretical point of view, developing a distributed setting to solve the tractable formulation \eqref{robust_comp_opt_multiagent}, e.g., extending the work in \cite{conte2013robust} to the  case  where binary variables are also present. We proposed a distributed stochastic MPC setting in \cite{rostampour2017distributedsmpc} and plan to extend such a framework to cope with a large-scale mixed-integer stochastic program, e.g., \eqref{robust_comp_opt_multiagent}.

c) Another interesting research direction is to integrate our developments in this paper into a power grid.
This may be a promising framework  such that ATES and   power grid  operators  can  simultaneously optimize  their performance to overcome uncertainties in generation \cite{rostampour2017tractable}, e.g., wind power, and  operational  demand decisions, the level of comfort for BCC systems.

\section*{Acknowledgment}
The authors are grateful to Wayan Wicak Ananduta for fruitful discussions and to Marc Jaxa-Rozen for his contributions to the simulation study depicted in Fig.~\ref{fig_marc}.

\appendix

\noindent\textbf{Proof of \Cref{thm_common}.}
Define $\text{Vio}^p(\bs{y}^*_{s})$, and $\text{Vio}^c(\bs{y}^*_{s})$ to be the violation probabilities of the private and common chance constraints as in \eqref{private_uncer}, and \eqref{common_uncer}, respectively, as follows:
%
\begin{align*}
\text{Vio}^p(\bs{y}^*_{s}) &= \Pb_{\bs{w}}  \left[ \bs{w} \in \cl{W} : \bs{y}^*_s \notin \prod_{i=1}^{N} \cl{Y}_{i}(\bs{w}_{i}) \ , \ \bs{y}^*_s \in \cl{Y} \right] \,, &
\\
\text{Vio}^c(\bs{y}^*_{s}) &= \Pb_{\bs{\delta}} \left[ \bs{\delta} \in \Delta : \bs{y}^*_s \notin \prod_{i=1}^{N} \bigcap_{j\in \cl{N}_i}^{} \breve{\cl{Y}}_{ij}(\bs{\delta}_{ij})  \ , \ \bs{y}^*_s \in \cl{Y} \right] \,, 
\end{align*}
%
where $\cl{Y}$ is the feasible region of the problem \eqref{robust_comp_opt_multiagent}, and it can be characterized via
\begin{align*}
\cl{Y} := \left\{ \bs{y}\in\R^{n_y N} \, : \, \bs{y}\in \left\{\prod\limits_{i=1}^{N} \bigcap\limits_{\bs{w}_i\in \{\cl{B}^*_i \bigcap \cl{W}_i\}} \cl{Y}_{i}(\bs{w}_{i}) \right\} 
\  \bigcap \ 
\left\{ \prod\limits_{i=1}^{N} \bigcap\limits_{j\in \cl{N}_i}^{} \bigcap\limits_{\bs{\delta}_{ij}\in \{\bar{\cl{B}}^*_{ij} \bigcap \Delta_{ij}\}} \, \breve{\cl{Y}}_{ij}(\bs{\delta}_{ij})  \right\}\,\right\} \,.
\end{align*}
It is important to note that both definitions of $\text{Vio}^p(\bs{y}^*_{s})$ and $\text{Vio}^c(\bs{y}^*_{s})$ have to be conditional probabilities $\Pb_{\bs{w}},  \left[(\cdot) | \forall\bs{\delta}\in\Delta \right]$, and $\Pb_{\bs{\delta}}  \left[(\cdot) | \forall\bs{w}\in\cl{W}\right]$, respectively.
We however equivalently considered them in the above form, due to the independency of both random process following Assumption \ref{independence}.
Define $\cl{B}^* = \prod\limits_{i=1}^{N} \cl{B}^*_i$, and consider now:
\begin{align*}
\bs{y}^* \in \prod_{i=1}^{N} \bigcap_{\bs{w}_i\in \{\cl{B}^*_i \bigcap \cl{W}_i\}} \cl{Y}_{i}(\bs{w}_{i}) 
\ \Leftrightarrow \
\bs{y}^* \in \bigcap_{\bs{w}\in \{\cl{B}^* \bigcap \cl{W}\}} \prod_{i=1}^{N} \cl{Y}_{i}(\bs{w}_{i})	\,.
\end{align*}
We also define $\bar{\cl{B}}^* = \prod\limits_{i=1}^{N} \bar{\cl{B}}^*_i$ and $\bar{\cl{B}}^*_i = \prod\limits_{j\in \cl{N}_i}^{} \bar{\cl{B}}^*_{ij}$, and clearly, we can have:
\begin{align*}
\bs{y}^* \in \prod\limits_{i=1}^{N} \bigcap\limits_{j\in \cl{N}_i}^{} \bigcap\limits_{\bs{\delta}_{ij}\in \{\bar{\cl{B}}^*_{ij} \bigcap \Delta_{ij}\}} \breve{\cl{Y}}_{ij}(\bs{\delta}_{ij})
\Leftrightarrow 
\bs{y}^* \in \prod\limits_{i=1}^{N} \bigcap\limits_{\bs{\delta}_{i}\in \{\bar{\cl{B}}^*_{i} \bigcap \Delta_{i}\}} \bigcap\limits_{j\in \cl{N}_i}^{} \breve{\cl{Y}}_{ij}(\bs{\delta}_{ij}) 
 \Leftrightarrow 
\bs{y}^* \in \bigcap\limits_{\bs{\delta}_{}\in \{\bar{\cl{B}}^*_{} \bigcap \Delta_{}\}} \prod\limits_{i=1}^{N} \bigcap\limits_{j\in \cl{N}_i}^{} \breve{\cl{Y}}_{ij}(\bs{\delta}_{ij}) \ .
\end{align*}
Therefore, if ${\bs{w}\in \{\cl{B}^* \bigcap \cl{W}\}}$ then $\bs{y}^* \in \prod\limits_{i=1}^{N} \, \cl{Y}_{i}(\bs{w}_{i})$, and 
if ${\bs{\delta}_{}\in \{\bar{\cl{B}}^*_{} \bigcap \Delta_{}\}}$ then $\bs{y}^* \in \prod\limits_{i=1}^{N} \, \bigcap\limits_{j\in \cl{N}_i}^{} \, \breve{\cl{Y}}_{ij}(\bs{\delta}_{ij})$.
This yields the following relations:
\begin{align*}
\text{Vio}^p(\bs{y}^*_{s}) &\leq \Pb_{\bs{w}}  \left[ \bs{w} \in \cl{W} : \bs{w}\notin \cl{B}^* \right] = \text{Vio}(\cl{B}^* ) \,, 
\\
\text{Vio}^c(\bs{y}^*_{s}) &\leq \Pb_{\bs{\delta}} \left[ \bs{\delta} \in \Delta : \bs{\delta}_{} \notin \bar{\cl{B}}^*_{} \right] = \text{Vio}(\bar{\cl{B}}^*_{}) \,, 
\end{align*}
It is then sufficient to show that for $N_s = \max_{i = 1, \cdots, N} N_{s_i}$, and $\bar{N}_s = \max_{i = 1, \cdots, N} \max_{j\in N_j} \bar{N}_{s_{ij}}$:
\begin{subequations}
	\label{desired_bound}
	\begin{align}
	\Pb^{N_s}_{\bs{w}} \Big[\cl{S} \in \cl{W}^{N_s} : \text{Vio}(\cl{B}^*) \geq \varepsilon \Big] &\leq  \beta \ , \label{desired_bound_w} \\
	\Pb^{\bar{N}_s}_{\bs{\delta}} \Big[\bar{\cl{S}} \in \Delta^{\bar{N}_s} : \text{Vio}(\bar{\cl{B}}^*_{}) \geq \bar\varepsilon \Big] &\leq  \bar\beta \ , \label{desired_bound_d}
	\end{align}
\end{subequations}
where $\cl{S} = \prod_{i=1}^{N} \cl{S}_i$, and $\bar{\cl{S}} = \prod_{i=1}^{N} \prod_{j\in \cl{N}_i}^{} \bar{\cl{S}}_{ij}$.
To this end, we now break down the proof in the following steps to show \eqref{desired_bound_w} and \eqref{desired_bound_d}: 
\begin{enumerate}
	
	\item Common chance constraint violation:
	\begin{align*}
	\text{Vio}^c(\bs{y}^*_{s}) \leq \text{Vio}(\bar{\cl{B}}^*_{})
	&= \Pb_{\bs{\delta}} \left[ \bs{\delta} \in \Delta : \bs{\delta}_{} \notin \bar{\cl{B}}^*_{} \right] \\
	&= \Pb_{\bs{\delta}}  \left[\bs{\delta} \in \Delta : \bs{\delta}_{} \notin \prod_{i=1}^{N} \prod_{j\in \cl{N}_i}^{} \bar{\cl{B}}^*_{ij} \right]   \\
	&= \Pb_{\bs{\delta}}  \Big[\bs{\delta} \in \Delta : \exists i \in\{1, \cdots, N\}, \bs{\delta}_{i} \notin  \prod_{j\in \cl{N}_i}^{} \bar{\cl{B}}^*_{ij} \Big] \\
	&= \Pb_{\bs{\delta}}  \left[\bigcup_{i=1}^{N}\left\{ \bs{\delta}_i \in \Delta_i : \bs{\delta}_{i} \notin \prod_{j\in \cl{N}_i}^{} \bar{\cl{B}}^*_{ij} \right\} \right]   \\
	&\leq \sum_{i=1}^{N} \Pb_{\bs{\delta}_i}  \left[ \bs{\delta}_i \in \Delta_i : \exists j\in \cl{N}_i \,,\,  \bs{\delta}_{ij} \notin \bar{\cl{B}}^*_{ij} \right] \\
	&= \sum_{i=1}^{N} \Pb_{\bs{\delta}_i}  \left[\bigcup_{j\in \cl{N}_i}^{}\left\{ \bs{\delta}_{ij} \in \Delta_{ij} : \bs{\delta}_{ij} \notin \bar{\cl{B}}^*_{ij}\right\} \right]   \\
	&\leq \sum_{i=1}^{N} \sum_{j\in \cl{N}_i} \Pb_{\bs{\delta}_{ij}} \left[\bs{\delta}_{ij} \in \Delta_{ij} : \bs{\delta}_{ij} \notin \bar{\cl{B}}^*_{ij} \right] \\
	&= \sum_{i=1}^{N} \sum_{j\in \cl{N}_i} \text{Vio}(\bar{\cl{B}}^*_{ij})\ .
	\end{align*}
	This implies that $\text{Vio}(\bs{y}^*_{s}) \leq \sum\limits_{i=1}^{N} \sum\limits_{j\in \cl{N}_i} \text{Vio}(\bar{\cl{B}}^*_{ij})$, and thus, we have
	\begin{align*}
	\Pb^{\bar{N}_s}_{\bs{\delta}}\left[\bar{\cl{S}} \in \Delta^{\bar{N}_s} : \text{Vio}(\bs{y}^*_{s}) \geq \bar\varepsilon\right]
	&\leq \Pb^{\bar{N}_s}_{\bs{\delta}} \left[\bar{\cl{S}} \in \Delta^{\bar{N}_s} : \sum\limits_{i=1}^{N} \sum\limits_{j\in \cl{N}_i} \text{Vio}(\bar{\cl{B}}^*_{ij}) \geq \sum_{i=1}^{N} \sum_{j\in \cl{N}_i}\bar\varepsilon_{ij} \right]  \\
	& = \Pb^{\bar{N}_s}_{\bs{\delta}} \left[ \bigcup_{i=1}^{N} \left\{\bar{\cl{S}}_i \in \Delta^{\bar{N}_{s_i}}_i :\sum_{j\in \cl{N}_i} \text{Vio}(\bar{\cl{B}}^*_{ij}) \geq \sum_{j\in \cl{N}_i}\bar\varepsilon_{ij} \right\} \right]   \\
	& \leq \sum_{i=1}^{N} \Pb^{\bar{N}_{s_i}}_{\bs{\delta}_i} \left[ \bar{\cl{S}}_i \in \Delta^{\bar{N}_{s_i}}_i :\sum_{j\in \cl{N}_i} \text{Vio}(\bar{\cl{B}}^*_{ij}) \geq \sum_{j\in \cl{N}_i}\bar\varepsilon_{ij} \right] \\  
	& = \sum_{i=1}^{N} \Pb^{\bar{N}_{s_i}}_{\bs{\delta}_i} \left[\bigcup_{j\in \cl{N}_i}^{} \left\{ \bar{\cl{S}}_{ij} \in \Delta^{\bar{N}_{s_{ij}}}_{ij} :\text{Vio}(\bar{\cl{B}}^*_{ij}) \geq \bar\varepsilon_{ij} \right\} \right] \\  
	& \leq \sum_{i=1}^{N} \sum_{j\in \cl{N}_i} \Pb^{\bar{N}_{s_{ij}}}_{\bs{\delta}_{ij}} \left[\bar{\cl{S}}_{ij} \in \Delta^{\bar{N}_{s_{ij}}}_{ij} :\text{Vio}(\bar{\cl{B}}^*_{ij}) \geq \bar\varepsilon_{ij} \right] \\
	&\leq \ \sum_{i=1}^{N} \sum_{j\in \cl{N}_i} \beta_{ij} = \beta \ . 
	\end{align*}

	\item Private chance constraint violation:
	\begin{align*}
	\text{Vio}^p(\bs{y}^*_{s}) \leq \text{Vio}(\cl{B}^*) 
	&= \Pb_{\bs{w}} \left[\bs{w} \in \cl{W} : \bs{w}\notin \cl{B}^* \right] \\
	&= \Pb_{\bs{w}} \left[\bs{w} \in \cl{W} : \bs{w} \notin \prod_{i=1}^{N} \cl{B}^*_i \right] \\
	&= \Pb_{\bs{w}} \left[\bs{w} \in \cl{W} : \exists i \in\{1, \cdots, N\} \, , \, \bs{w}_i \notin \cl{B}^*_i \right] \\
	&= \Pb_{\bs{w}} \left[\bigcup_{i=1}^{N}\left\{ \bs{w}_i \in \cl{W}_i : \bs{w}_i \notin \cl{B}^*_i \right\} \right] \\
	&\leq \sum_{i=1}^{N} \Pb_{\bs{w}_i} \left[ \bs{w}_i \in \cl{W}_i : \bs{w}_i \notin \cl{B}^*_i \right] \\
	&= \sum_{i=1}^{N} \text{Vio}(\cl{B}^*_i) \ . 
	\end{align*}
	The last statement implies that $\text{Vio}^p(\bs{y}^*_{s}) \leq \sum\limits_{i=1}^{N} \text{Vio}(\cl{B}^*_i)$, and thus, we have
	\begin{align*}
	\Pb^{N_s}_{\bs{w}} \left[\cl{S} \in \cl{W}^{N_s} : \text{Vio}^p(\bs{y}^*_{s}) \geq \varepsilon\right]
	&\leq  
	\Pb^{N_s}_{\bs{w}} \Big[\cl{S} \in \cl{W}^{N_s} : \sum_{i=1}^{N} \text{Vio}(\cl{B}^*_i) \geq \sum_{i=1}^{N}\varepsilon_i \Big]  \\
	& = \Pb^{N_{s}}_{\bs{w}} \left[ \bigcup_{i=1}^{N} \left\{ \cl{S}_i \in \cl{W}^{N_{s_i}}_i : \text{Vio}(\cl{B}^*_i) \geq \varepsilon_i \right\} \right]   \\
	& \leq \sum_{i=1}^{N} \Pb^{N_{s_i}}_{\bs{w}_i} \left[ \cl{S}_i \in \cl{W}^{N_{s_i}}_i : \text{Vio}(\cl{B}^*_i) \geq \varepsilon_i \right] \\
	&\leq  \ \sum_{i=1}^{N} \beta_i = \beta \ . 
	\end{align*}

\end{enumerate}
The obtained bounds in the above procedure are the desired assertions as it is stated in the theorem.
It is important to mention that we use the existing results in \cite{campi2008exact} to determine ${N}_{s_{i}}$ and $\bar{N}_{s_{ij}}$ and solve the tractable problems \eqref{SPBi} and \eqref{SPBij} for each agent $i = 1, \cdots, N$, $\forall j \in \cl{N}_i$, respectively. 
We thus have the following probabilistic guarantees: 
\begin{align*}
&\Pb^{N_{s_i}}_{\bs{w}_i} \left[ \cl{S}_i \in \cl{W}^{N_{s_i}}_i : \text{Vio}(\cl{B}^*_i) \geq \varepsilon_i \right] \ \leq \beta_i \ , \\
&\Pb^{\bar{N}_{s_{ij}}}_{\bs{\delta}_{ij}} \left[\bar{\cl{S}}_{ij} \in \Delta^{\bar{N}_{s_{ij}}}_{ij} :\text{Vio}(\bar{\cl{B}}^*_{ij}) \geq \bar\varepsilon_{ij} \right] \ \leq \beta_{ij} \ .
\end{align*}
The interpretation of the derivation of these bounds \eqref{desired_bound} is as follows.
The probability of all violation probabilities $\text{Vio}(\cl{B}^*_i)$ being simultaneously bounded by the corresponding $\varepsilon_i$ is at least $1-\beta$,  
and $\text{Vio}(\bar{\cl{B}}^*_{ij})$ being simultaneously bounded by the corresponding $\bar\varepsilon_{ij}$ is at least $1-\bar\beta$. 
The proof is completed by noting that the feasible set  $\cl{Y}$ of \eqref{robust_comp_opt_multiagent} has a non-empty interior:
\begin{align*}
\Big\{\exists \rho \in \R_{+} \ , \ \bar{\bs{y}} \in \cl{Y}  \ : \ \|\bs{y} - \bar{\bs{y}}\| \leq \rho , \, \forall \bs{y}\in\R^{n_y}\Big\} \subset \cl{Y} \ ,
\end{align*}
and since the problem \eqref{robust_comp_opt_multiagent} has a non-empty interior feasible set, it admits at least one feasible solution $\bs{y}^*_{s}$.
\QEDB


\noindent\textbf{Proof of \Cref{cor_common_cons}.}
The proof is straightforward by just substituting the corresponding relationships, we have 
\begin{align*}
(\tx{V}_{\tx{a},k}^h)_i + (\tx{V}_{\tx{a},k}^c)_j &\leq \tx{V}_{ij} - \bar\delta_{ij,k} \,, \\
\frac{c_{aq}  \pi \ell}{c_{pw}} \Big((\bar{r}_{\tx{a},k}^h)_i^2 + (\bar{r}_{\tx{a},k}^c)_j^2 \Big) &\leq \frac{c_{aq}  \pi \ell}{c_{pw}} \Big( (d_{ij})^2 - 2 (\bar{r}_{\tx{a},k}^h)_i \, (\bar{r}_{\tx{a},k}^c)_j \Big)\,,  \\
(\bar{r}_{\tx{a},k}^h)_i^2 + (\bar{r}_{\tx{a},k}^c)_j^2  &\leq (d_{ij})^2 - 2 (\bar{r}_{\tx{a},k}^h)_i \, (\bar{r}_{\tx{a},k}^c)_j \,, \\
\Big( (\bar{r}_{\tx{a},k}^h)_i + (\bar{r}_{\tx{a},k}^c)_j \Big)^2  &\leq (d_{ij})^2  \,,  \\
(\bar{r}_{\tx{a},k}^h)_i + (\bar{r}_{\tx{a},k}^c)_j &\leq d_{ij}\,.
\end{align*}
The proof is completed by noting that the thermal radius is positive: $(\bar{r}_{\tx{a},k}^h)_i \geq 0\,,$ $(\bar{r}_{\tx{a},k}^c)_i \geq 0$, $\forall i \in\{ 1, \cdots, N\}$.
\QEDB

\bibliographystyle{IEEEtran}
\bibliography{references}

\end{document}